\documentclass[11pt]{amsart}
\usepackage{amsmath,amsfonts,latexsym,amssymb,amsthm,mathrsfs,array}

\usepackage{amsfonts}
\usepackage[usenames]{color}
\usepackage[colorlinks=true,linkcolor=webgreen,filecolor=webbrown,citecolor=webgreen]
{hyperref}
\definecolor{webgreen}{rgb}{0,.5,0}
\definecolor{webbrown}{rgb}{.6,0,0}

\newcommand{\CC}{{\mathbb C}}
\newcommand{\ZZ}{{\mathbb Z}}

\def\gg{{\mathfrak{g}}}
\def\hh{{\mathfrak{h}}}
\def\CC{{\mathbb C}}
\def\ZZ{{\mathbb Z}}
\newcommand{\fma}{\overset{\circ}{\mathfrak{g}}}

\newtheorem{dfn}{Definition}[section]
\newcommand{\bdfn}{\begin{dfn}\rm}
\newcommand{\edfn}{\end{dfn}}
\newtheorem{thm}[dfn]{Theorem}
\newcommand{\bthm}{\begin{thm}}
\newcommand{\ethm}{\end{thm}}
\newtheorem{lmma}[dfn]{Lemma}                   
\newcommand{\blmma}{\begin{lmma}}                   
\newcommand{\elmma}{\end{lmma}}                   
\newtheorem{ppsn}[dfn]{Proposition}
\newcommand{\bppsn}{\begin{ppsn}}
\newcommand{\eppsn}{\end{ppsn}}
\newtheorem{crlre}[dfn]{Corollary}
\newcommand{\bcrlre}{\begin{crlre}} 
\newcommand{\ecrlre}{\end{crlre}}
\newtheorem{rmk}[dfn]{Remark}
\newcommand{\brmk}{\begin{rmk}\rm} 
\newcommand{\ermk}{\end{rmk}}

\numberwithin{equation}{section}

\title[Integrable Modules for twisted full toroidal Lie algebras] {Classification of level zero Irreducible Integrable modules for twisted full toroidal Lie algebras}
\author{Souvik Pal}
\author{S. Eswara Rao}        
\address{Souvik Pal, Harish-Chandra Research Institute (HBNI), Chhatnag Road, Jhunsi, Prayagraj(Allahabad) 211019, Uttar Pradesh, India.}     
\email{pal.souvik90@gmail.com, souvikpal@hri.res.in}     
\address{S. Eswara Rao, School of Mathematics, Tata Institute of Fundamental Research,
Homi Bhabha Road, Colaba, Mumbai 400005, India.}                  
\email{sena98672@gmail.com, senapati@math.tifr.res.in}   
  
\date{}

\begin{document}

\maketitle
\begin{abstract}
In this paper, we first construct the twisted full toroidal Lie algebra by 
an extension of a centreless Lie torus $LT$ which is a multiloop algebra twisted by several automorphisms of finite order and equipped with a particular grading. We then provide a complete classification of all the irreducible integrable modules with finite-dimensional weight spaces for this twisted full toroidal Lie algebra having a non-trivial $LT$-action and where the centre of the underlying Lie algebra acts trivially.\\\\                                                                    {\bf{MSC}:} Primary: 17B67; Secondary: 17B65, 17B70.\\
{\bf{KEY WORDS}:} Lie torus, integrable, twisted toroidal, level zero.           
\end{abstract} 

\section{Introduction}

It is well-known that the Virasoro algebra plays a prominent role in the representation theory of affine Lie algebras as it acts on $almost$ every highest weight module for the affine Lie algebra via the famous Sugawara operators. This remarkable connection eventually resulted in the construction of the so-called affine-Virasoro algebra \cite{V,G} which is the semi-direct product of the Virasoro algebra and the derived affine Kac-Moody algebra.   Subsequently this Lie algebra has emerged to be an extremely important object of study in various branches of mathematics and physics. In particular, its relation to Conformal Field Theory has been explained in great detail in \cite{FMS}. The twisted multivariable generalization of this classical object is the subject of our current paper.

Our pursuit of a suitable replacement for the twisted affine Lie algebra leads us to the notion of a $centreless$ $Lie$ $torus$ satisfying $fgc$ $condition$ (which means that the Lie torus is finitely generated as a module over its centroid). It was shown in \cite{ABFP} that such a Lie torus of non-zero nullity has a $multiloop$ $realization$. Moreover in the one variable case, this multiloop algebra is simply the loop algebra twisted by a single Dynkin automorphism used in the construction of twisted affine Kac-Moody algebras \cite{K}. It is worth mentioning here that centreless Lie tori satisfying fgc condition play a key role in the theory of extended affine Lie algebras (EALAs for short) as they turn out to be the $centreless$ $cores$ $of$ $almost$ $every$ EALA \cite{Ne}.

We first start with a finite-dimensional simple Lie algebra $\mathfrak{g}$ endowed with a specific grading (see (\ref{Eigenspace})). We then introduce the centreless Lie torus $LT$ (see Definition \ref{Tori}) and henceforth consider its universal central extension $\overline{LT}= LT \oplus \mathcal{Z}(\underline{m})$ (see Subsection \ref{Extension}). The next step is to generalize the Virasoro algebra in the multivariable setup. To this end, let us consider the algebra of derivations on $A(\underline{m}) = \CC[t_1^{\pm m_1}, \ldots, t_n^{\pm m_n}]$ which we shall denote by $\mathcal{D}(\underline{m})$ (see (\ref{Witt})). Unlike in the single variable case, this derivation algebra is centrally closed for $n \geqslant 2$ \cite{RSS}. Nonetheless, it admits an abelian extension given by the centre $\mathcal{Z}(\underline{m})$ of $\overline{LT}$ which generalizes the Virasoro algebra \cite{EM}. We can now take the semi-direct product of $\overline{LT}$ with $\mathcal{D}(\underline{m})$ and finally obtain the twisted full toroidal Lie algebra $\widetilde{\tau}$ (see Subsection \ref{Twisted Toroidal}). The centre of this resulting Lie algebra is spanned by only finitely many elements $K_1, \ldots, K_n$. If all these $n$ elements act trivially on a $\widetilde{\tau}$-module, we say that the $level$ of this representation is $zero$, otherwise we say that the representation has $non$-$zero$ $level$.

The classification of irreducible integrable modules with finite-dimensional weight spaces for the untwisted as well as for the twisted affine Kac-Moody algebras have been carried out in \cite{C,CP2,CP1,CP3}. In the multivariable setup, if the $n$-tuple of automorphisms $\sigma=(\sigma_1, \ldots, \sigma_n)$ corresponding to our Lie torus $L(\mathfrak{g}, \sigma)$ (see Definition \ref{D2.1}) are all taken to be identity, then $\widetilde{\tau}$ is simply the $untwisted$ full toroidal Lie algebra (FTLA for short). In this context, a large class of modules has been explicitly constructed in \cite{B,YB} through the use of vertex operator algebras. The irreducible integrable modules with finite-dimensional weight spaces for FTLA have been classified in \cite{EJ}.

The classification of some $specific$ irreducible integrable representations of $non$-$zero$ $level$ with finite-dimensional weight spaces for the multiloop algebra twisted by a single non-trivial automorphism first appeared in \cite{FJ}. Thereafter this was generalized for all irreducible integrable modules of $non$-$zero$ $level$ over the twisted full toroidal Lie algebra (TFTLA for short) in \cite{BE}. Realizations of several modules of $non$-$zero$ $level$ have also been provided in \cite{BL} for the twisted case. The current paper addresses these problems for the first time in case of $level$ $zero$ modules for TFTLA.

The main purpose of our paper is to classify the irreducible integrable level zero modules with finite-dimensional weight spaces for TFTLA having a non-trivial action of the Lie torus. If the action of $LT$ is trivial, then these modules eventually turn out to be irreducible over $\mathcal{D}(\underline{m})$ which have been classified in \cite{BF}. We now summarize the content of our paper.

After constructing the TFTLA in Section \ref{Notation}, we observe that this  Lie algebra $\widetilde{\tau}$ has a natural triangular decomposition (see (\ref{LT})) given by
\begin{align*}
\widetilde{\tau} = \widetilde{\tau}_{-} \oplus \widetilde{\tau}_{0} \oplus \widetilde{\tau}_{+}.         
\end{align*}
Suppose that $V$ is an irreducible integrable level zero module over $\widetilde{\tau}$ with finite-dimensional weight spaces. 
In Section \ref{Representations}, we prove that $V$ is actually a highest weight module in the sense of Definition \ref{D3.3} (see Proposition \ref{Highest}). In Section \ref{Central}, we show that the action of $\mathcal{Z}(\underline{m})$ is trivial on our module $V$ (see Corollary \ref{C4.9}). We then recall some results from \cite{EK} and \cite{L} in the following section and thereby utilize them in Section \ref{Realization} where we construct a family of level zero irreducible integrable modules with finite-dimensional weight spaces for TFTLA and having a non-trivial $LT$-action (see Proposition \ref{P6.1}). These are the first such realizations of level zero irreducible modules for TFTLA. In Section \ref{Classification}, we show that this collection completely exhausts all such irreducible integrable modules over TFTLA (see Corollary \ref{C7.11}). We then proceed to Section \ref{Class} where we find the isomorphism classes of the aforementioned irreducible integrable modules (see Theorem \ref{T8.4}). In the last section, we obtain the main result of our paper which conveys that the level zero irreducible integrable modules with finite-dimensional weight spaces and having a non-trivial action of the Lie torus are completely as well as uniquely determined by a finite-dimensional irreducible $\mathfrak{gl}_n(\mathbb{C})$-module and a finite-dimensional graded-irreducible $\mathfrak{g}$-module up to a shift of grading.                                                                                            

We would like to emphasize that many of our results actually hold good for a more general class of irreducible highest weight modules that need not be integrable (for instance, see Theorem \ref{T4.1}). These can be used in the study of highest weight modules that may possibly help us to understand the famous $Harish$-$Chandra$ $modules$. The irreducible Harish-Chandra modules have been classified over some well-known Lie algebras \cite{FT,ML,M,Ma}. Finally observe that our TFTLA is intimately connected to the twisted toroidal EALA (TTEALA for short) which is formed by just adjoining the subalgebra $\mathcal{S}(\underline{m})$= span$\{D(\underline{u}, \underline{r}) := \sum_{i=1}^{n}u_it^{\underline{r}}d_i \ | \ (\underline{u}, \underline{r})=0, \ \underline{u} \in \mathbb{C}^n, \ \underline{r} \in \Gamma \}$ of $\mathcal{D}(\underline{m})$ to $\overline{LT}$, where $(\cdot, \cdot)$ denotes the standard inner product on $\mathbb{C}^n$ (see \cite{BES} for more details on TTEALAs). The representation theory of TTEALAs is still in progress and most of the work has been done only with regard to representations of non-zero
level \cite{BL,BES}. So we hope that this paper will also contribute to develop the representation theory of EALAs, especially in case of level zero modules.              
\section{Notations and Preliminaries}\label{Notation}
Throughout the paper, all the vector spaces, algebras and tensor products are over the field of complex numbers $\mathbb{C}$. We shall denote the set of integers, natural numbers, non-negative integers and non-zero complex numbers by $\mathbb{Z}$, $\mathbb{N}$, $\mathbb{Z_{+}}$ and $\mathbb{C}^{\times}$ respectively. For any Lie algebra $L$, its universal enveloping algebra will be denoted by $U(L)$.  
\subsection{Full Toroidal Lie Algebras} \label{Toroidal}
Consider a finite-dimensional simple Lie algebra $\mathfrak{g}$ equipped with a Cartan subalgebra $\mathfrak{h}$. Then $\mathfrak{g}$ is endowed with a symmetric, non-degenerate bilinear form which remains invariant under every automorphism of $\mathfrak{g}$. We shall denote this bilinear form by $(\cdot|\cdot)$.\\ 
Let $A = \CC[t_1^{\pm 1}, \ldots, t_n^{\pm 1}]$ be the algebra of Laurent polynomials in $n$ variables.
Consider the (untwisted) multiloop algebra given by
\begin{align*}
L(\mathfrak{g})=\mathfrak{g} \otimes A 
\end{align*}
which is a Lie algebra under the following pointwise bilinear operation.
\begin{align*}
[x \otimes f, y \otimes g] = [x,y] \otimes fg \ \forall \ x,y \in \mathfrak{g} \ \text{and} \ f, g \in A.
\end{align*}
For any $x \in \mathfrak{g}$ and $\underline{k} \in \mathbb{Z}^n$, write $t^{\underline{k}}=t_1^{k_1} \ldots t_n^{k_n}$ and let $x(\underline{k}) = x \otimes t^{\underline{k}}$ denote a typical element of $L(\mathfrak{g})$. Now consider the module of differentials ($\Omega_A,d$) of $A$ which is the free $A$-module with basis $\{K_1, \ldots, K_n \}$ along with the differential map $d : A \longrightarrow \Omega_A$ whose image is the subspace spanned by the elements of the form $\sum_{i =1}^{n}k_i t^{\underline{k}}K_i$ for $\underline{k} \in \ZZ^{n}$. More precisely, we have       
\begin{align*}
	\Omega_A = \text{span} \{t^{\underline{k}}K_i \ | \ 1 \leq i \leq n, \ \underline{k} \in \ZZ^{n} \},\ 
	dA =  \text{span} \{\sum_{i =1}^{n}{k_i t^{\underline{k}}K_i} \ | \ \underline{k} \in \ZZ^{n} \}.
\end{align*}
If we now consider the quotient space $\mathcal{Z} = \Omega_A / dA$, then we know that
\begin{align}\label{Universal}
	\overline{L}(\mathfrak{g})=L(\mathfrak{g}) \oplus \mathcal{Z}
\end{align} 
is the universal central extension of $L(\mathfrak{g})$ \cite{Ka,EMY}. By abuse of notation, we shall denote the image of $t^{\underline{k}}K_i$ in $\mathcal{Z}$ again by itself. Henceforth $\overline{L}(\mathfrak{g})$ also forms a Lie algebra under the following bracket operations. 
\begin{enumerate}
	\item $[x(\underline{k}), y(\underline{l})] = [x,y](\underline{k} + \underline{l}) + (x|y) \sum_{i=1}^{n} k_it^{\underline{k} + \underline{l}}K_i$,
	\item $\mathcal{Z}$ is central in $\overline{L}(\mathfrak{g})$.  
\end{enumerate}
The final ingredient needed to define the (untwisted) full toroidal Lie algebra is the Lie algebra of derivations on $A$ which we shall denote by $\mathcal{D}$. Taking $d_i=t_i \frac{d}{dt_i}$ for $1 \leqslant i \leqslant n$ which acts on $A$ as derivations,  we can also define
\begin{align*}
\mathcal{D}= \text{span} \{t^{\underline{r}}d_i : \underline{r} \in \mathbb{Z}^n, \ 1 \leqslant i \leqslant n \}.    
\end{align*}
This infinite-dimensional derivation algebra is an extremely classical object in its own right and is popularly known as the multivariable $Witt$ $algebra$. It is now easy to verify that
\begin{align*}
[t^{\underline{r}}d_i, t^{\underline{s}}d_j]=s_it^{\underline{r} + \underline{s}}d_j - r_jt^{\underline{r} + \underline{s}}d_i.
\end{align*}
Again any $d \in \mathcal{D}$ can be extended to a derivation on $L(\mathfrak{g})$ by simply setting
\begin{align*}
d(x \otimes f)=x \otimes df \ \forall \ x \in \mathfrak{g}, \ f \in A
\end{align*}
which subsequently has a unique extension to $\overline{L}(\mathfrak{g})$ via the
action  
\begin{align*}
t^{\underline{r}}d_i(t^{\underline{s}}K_j)=s_it^{\underline{r} + \underline{s}}K_j + \delta_{ij} \sum_{p=1}^{n}r_p t^{\underline{r} + \underline{s}}K_p.
\end{align*}
Moreover it is well-known that $\mathcal{D}$ admits two non-trivial 2-cocyles $\phi_1$ and $\phi_2$ with values in $\mathcal{Z}$ (see \cite{BB} for more details) where we have
\begin{align*}
\phi_1(t^{\underline{r}}d_i, t^{\underline{s}}d_j) = -s_ir_j \sum_{p=1}^{n}{r_p t^{\underline{r} + \underline{s}}K_p}, 
\end{align*}
\begin{align*}
\phi_2(t^{\underline{r}}d_i, t^{\underline{s}}d_j) = r_is_j \sum_{p=1}^{n}{r_p t^{\underline{r} + \underline{s}}K_p}.          
\end{align*}
Let $\phi$ be any arbitrary linear combination of $\phi_1$ and $\phi_2$. Then we can define the (untwisted) full toroidal Lie algebra (relative to $\mathfrak{g}$ and $\phi$) by setting
\begin{align}
\tau_A= L(\mathfrak{g}) \oplus \mathcal{Z} \oplus \mathcal{D}
\end{align}
with the following bracket operations besides the commutator relations (1) and (2) recorded in Subsection \ref{Toroidal}.
\begin{enumerate}
\item $[t^{\underline{r}}d_i,t^{\underline{s}}K_j]=s_it^{\underline{r} + \underline{s}}K_j + \delta_{ij} \sum_{p=1}^{n}r_p t^{\underline{r} + \underline{s}}K_p$,
\item $[t^{\underline{r}}d_i,t^{\underline{s}}d_j]=s_it^{\underline{r} + \underline{s}}d_j-r_j t^{\underline{r} + \underline{s}}d_i + \phi(t^{\underline{r}}d_i,t^{\underline{s}}d_j)$,
\item $[t^{\underline{r}}d_i,x \otimes t^{\underline{s}}]=s_ix \otimes t^{\underline{r} + \underline{s}} \ \forall \ x \in \mathfrak{g},\ \underline{r},\ \underline{s} \in \mathbb{Z}^n, \ 1 \leqslant i,j \leqslant n$.   
\end{enumerate}
\subsection{Lie Tori and Twisted Full Toroidal Lie Algebras}\label{Twisted}
Let us fix any $n\in\mathbb{N}$ and suppose that we have $n$ commuting automorphisms of a finite-dimensional simple Lie algebra $\mathfrak{g}$  given by $\sigma_1, \ldots, \sigma _n$ with finite orders $m_1, \ldots, m_n$ respectively. Put
\begin{align*} 
\sigma = (\sigma_1, \ldots, \sigma_n),\,\,\,\,\,\,\,\
\Gamma = m_1 \mathbb{Z} \oplus \ldots \oplus m_n\mathbb{Z},\,\,\,\,\,\,\,\
 G = \mathbb{Z}^n/ \Gamma.
 \end{align*} 
Thus we have a natural map \ $\mathbb{Z}^n \longrightarrow G$ $(\cong \mathbb{Z}/m_1\mathbb{Z} \times \ldots \times\mathbb{Z}/m_n\mathbb{Z})$ 
\begin{align}\label{Natural}
(k_1, \ldots, k_n) = \underline{k} \mapsto \overline{k}=(\overline{k_1}, \ldots ,\overline{k_n})
\end{align} 
For $1 \leqslant i \leqslant n$, let  $\xi_i$ denote a $m_i-$th primitive root of unity. Then we obtain an eigenspace decomposition of $\mathfrak{g}$ given by
\begin{align}\label{Eigenspace}
\mathfrak{g} = \bigoplus _{\overline{k} \in G} \mathfrak{g}_{\overline{k}} \,\,\,\ \text{where} \,\,\,\ \mathfrak{g}_{\overline{k}} = \{x \in \mathfrak{g} \ | \ \sigma_i x = \xi_i^{k_i}x,\ 1 \leqslant i \leqslant n \}.
\end{align}
It is well-known that $\mathfrak{g}_{\overline{0}}$ is a reductive Lie algebra (see \cite[Proposition 4.1]{BM}), even with the possibility of being zero.  Finally, let us define
\begin{align*} 
L(\mathfrak{g}, \sigma) = \bigoplus_{\underline{k} \in \mathbb{Z}^n} \mathfrak{g}_{\overline{k}} \otimes \mathbb{C}t^{\underline{k}}
\end{align*}
which is clearly a Lie subalgebra of the (untwisted) multiloop algebra $L(\mathfrak{g})$.  
\begin{dfn}\label{D2.1}
The above Lie algebra $L(\mathfrak{g}, \sigma)$ is known as the (twisted) multiloop algebra associated to $\mathfrak{g}$ and $ \sigma$.
\end{dfn}
For any finite-dimensional simple Lie algebra $\gg_1$ with a Cartan subalgebra $\hh_1$ of $\gg_1$, we have the root space decomposition of $\gg_1$ relative to $\hh_1$ given by
\begin{align*}
\gg_1 = \bigoplus_{\alpha \in \hh_1^*}\mathfrak{g}_{1, \alpha} \,\,\,\  \text{where} \,\,\,\  \mathfrak{g}_{1, \alpha} = \{ x \in \gg_1 \ | \ [h,x]= \alpha(h)x \ \forall \ h \in \hh_1\}.     
\end{align*}
Set $\Delta(\gg_1, \hh_1) = \{ \alpha \in \hh_1^*\ |\ \mathfrak{g}_{1, \alpha} \neq (0)\}$. Then $\Delta_{1}^{\times} = \Delta(\gg_1, \hh_1) \backslash \{0\}$ is clearly an irreducible reduced finite
root system having at most two root lengths. Let ${\Delta_{1, \mathrm{sh}}^{\times}}$ denote the set of all non-zero short roots of $\Delta_{1}^{\times}$. Define

\[
    {\Delta_{1, \mathrm{en}}^{\times}}= 
\begin{cases}
    \Delta_{1}^{\times} \cup 2  {\Delta_{1, \mathrm{sh}}^{\times}},\,\,\,\,\,\,\,\,\,\,\,\,\,\ &\text{if} \,\,\,\,\Delta_{1}^{\times} \ \text{is \ of \ type} \ B_l \\
    \Delta_{1}^{\times}  \,\,\,\,\,\,\,\,\,\,\,\,\,\,\,\,\,\,\,\,\,\ ,           & \text{otherwise}.   
\end{cases}
\]\
Suppose that $\{\alpha_1, \ldots ,\alpha_p\}$ is a collection of simple roots of $\mathfrak{g}_1$ with respect to $\mathfrak{h}_1$.         
Then $Q_1 = \bigoplus_{i=1}^{p} \mathbb{Z} \alpha_i$ is the corresponding root lattice containing the non-negative root lattice $Q_1^+ = \bigoplus_{i=1}^{p} \mathbb{Z}_{+} \alpha_i$. Also let  ${\Delta_{1, \mathrm{en}}^{\times,+}}$ and  ${\Delta_{1, \mathrm{en}}^{\times,-}}$ denote the respective set of positive and negative roots of the extended root system. Also set $\Delta_{1, \mathrm{en}} = {\Delta_{1, \mathrm{en}}^{\times}} \cup \{0\}$.\\    
For $\lambda, \mu \in \mathfrak{h}_1^*$, define a partial order relation $"\leqslant_1"$ on $\mathfrak{h}_1^*$ by setting\\
$\lambda \leqslant_1 \mu\ \text{if\ and\ only\ if}\ \mu-\lambda= \sum_{i=1}^{p} n_i \alpha_i\ \text{for\ some}\ n_1, \ldots ,n_p\ \text{in}\ \mathbb{Z}_{+}$.\\
Moreover, if $(\cdot,\cdot)$ is a non-degenerate, symmetric and associative bilinear form on $\mathfrak{g}_1$, then we shall denote the set of all dominant integral weights of $\mathfrak{g}_1$ with respect to $\mathfrak{h}_1$ by
\begin{align*}
P^{+}_{{\mathfrak{g}}_1} = \{ \lambda\ \in \mathfrak{h}_1^*\ | \ (\lambda,\alpha_i) \in \mathbb{Z}_{+}\ \forall \ i=1, \ldots ,p \}.
\end{align*}

Let us now define the notion of a Lie torus, which plays a central role in the construction of twisted full toroidal Lie algebras. In this paper, we shall use the characterization of a Lie torus given in \cite{ABFP} and take it as our definition rather than the axiomatic definition recorded in \cite{Y1,Y2}. For this purpose, we first need the following definition.
\bdfn\cite[Definition 3.2.3]{ABFP}
A finite-dimensional $\gg_1$-module $V$ is said to satisfy condition $(M)$ if 
\begin{enumerate}
\item $V$ is an irreducible module having dimension greater than one; 
\item The weights of $V$ relative to $\hh_1$ are contained in $ { \Delta_{1, \mathrm{en}}}$.
\end{enumerate}
\edfn

\bdfn\label{Tori}\cite[Proposition 3.2.5]{ABFP}
A multiloop algebra $L(\mathfrak{g}, \sigma)$ is called a Lie torus
if
\begin{enumerate}
\item $\mathfrak{g}_{\bar{0}}$ is a finite-dimensional simple Lie algebra;
\item For $\bar{k} \neq \overline{0}$ and $\mathfrak{g}_{\bar{k}} \neq (0)$, $\mathfrak{g}_{\bar{k}} \cong U_{\bar{k}} \oplus W_{\bar{k}}$, where $U_{\bar{k}}$
is trivial as a $\gg_{\bar{0}}$-module and either $W_{\bar{k}}$ is zero or it satisfies condition $(M)$;
\item $|<\sigma_1, \ldots , \sigma_n>|= \prod_{i=1}^{n} |\sigma_i|$  where $|\sigma_i|$ denotes the order of the automorphism $\sigma_i\ (1 \leqslant i \leqslant n)$ and $|<\sigma_1, \ldots , \sigma_n>|$ is the order of the group generated by the $\sigma_i$'s.
\end{enumerate}
\edfn

For the rest of this paper, we shall denote the Lie torus $L(\mathfrak{g}, \sigma)$ by $LT$. Now by \cite[Lemma A.1]{BL}, we can choose a Cartan subalgebra $\mathfrak{h}_{\bar{0}}$ of $\mathfrak{g}_{\overline{0}}$ such that $\mathfrak{h}_{\bar{0}} \subseteq 
\mathfrak{h}$. Then $\mathfrak{h}_{\bar{0}}$ turns out to be an ad-diagonalizable subalgebra of $\mathfrak{g}$ (see \cite[Lemma 3.1.3] {N}) and $\Delta^{\times} = \Delta^{\times}(\mathfrak{g}, \mathfrak{h}_{\bar{0}})$ is an irreducible (possibly non-reduced) finite root system (see \cite[Proposition 3.3.5]{N}). Let $Q$ be the root lattice corresponding to $\Delta^{\times}$. Finally put $\Delta_{0} = \Delta(\mathfrak{g}_{\bar{0}},\mathfrak{h}_{\bar{0}})$ and $\Delta = \Delta^{\times} \cup \{0\}$.          

\subsection{Properties of LT}\label{Torus} \cite{ABFP}
\begin{itemize}
\item[(LT1)] The Lie torus $LT$ is $\mathbb{Z}^n$-graded as well as $Q$-graded. This is usually referred to as a Lie $\mathbb{Z}^n$-torus of type $\Delta$.
\item[(LT2)]\label{Torus2}
\[\ {\Delta}= 
\begin{cases}
    \Delta_{0, \mathrm{en}}\ , &\text{if} \,\,\,\,\Delta_{0}^{\times} \ \text{is \ of \ type} \ B_l \\
    \Delta_0  \,\,\,\,\,\,\ ,           & \text{otherwise}.      
\end{cases}
\]
 \item[(LT3)] For each $\overline{k} \in G$, we have 
\begin{align*}
\mathfrak{g}_{\overline{k}} = \bigoplus_{\alpha \in \mathfrak{h}_{\overline{0}}^*} \mathfrak{g}_{\overline{k}}(\alpha)
\end{align*}
where $\mathfrak{g}_{\overline{k}}(\alpha) = \{ x \in \mathfrak{g}_{\overline{k}} \ | \ [h,x] = \alpha(h)x \ \forall\  h \in \mathfrak{h}_{\overline{0}}\}$ and dim$\mathfrak{g}_{\overline{k}}(\alpha) \leqslant1$ for every non-zero $\alpha$ in $\mathfrak{h}_{\overline{0}}^*$.

\item[(LT4)] The centre of $LT$ is trivial.
\end{itemize}

\subsection{Universal Central Extension of Lie Torus}\label{Extension}
Let $LT$ be the Lie $\mathbb{Z}^n$-torus of type $\Delta$ as defined earlier. Now just as in the case of $\tau_{A}$ (see Subsection \ref{Toroidal}), we can similarly construct the spaces $\Omega_{A(\underline{m})}$ and $dA(\underline{m})$ for the smaller algebra $A(\underline{m})= \mathbb{C}[t_1^{\pm m_1}, \ldots, t_n^{\pm m_n}]$. Set $\mathcal{Z}(\underline{m})= \Omega_{A(\underline{m})}/dA(\underline{m})$ and define
\begin{align*}
\overline{LT} := LT \oplus \mathcal{Z}(\underline{m}).  
\end{align*}
Then $\overline{LT}$ forms a Lie algebra under the following bracket operations.
\begin{enumerate}
\item $[x(\underline{k}), y(\underline{l})] = [x,y](\underline{k} + \underline{l}) \ + \ (x|y) \displaystyle{\sum_{i =1}^{n}} k_it^{\underline{k} + \underline{l}}K_i$ ;
\item $\mathcal{Z}(\underline{m})$ is central in $\overline{LT}$.
\end{enumerate}
\brmk\label{Generate}
\
\begin{enumerate}
\item Note that since $(\cdot|\cdot)$ is invariant under automorphisms of $\mathfrak{g}$, it is evident that $\underline{k} + \underline{l} \in \Gamma$ whenever $(x|y) \neq 0$ and thus the above bracket operation on $\overline{LT}$ is well-defined. 
\item $\overline{LT}$ is the universal central extension of $LT$ (see \cite[Corollary 3.27]{J}).    
\item $\overline{LT}$ is also a Lie $\mathbb{Z}^n$-torus of type $\Delta$ (with respect to the axiomatic definition of a Lie torus recorded in \cite[Definition 4.2]{EN}). This readily implies that $\overline{LT}$ is generated as a Lie algebra by $\bigcup_{\alpha \in \Delta^{\times}}(\overline{LT})_{\alpha}$, where $(\overline{LT})_{\alpha} = \bigoplus_{\underline{k} \in \mathbb{Z}^n} \mathfrak{g}_{\overline{k}}(\alpha) \otimes \mathbb{C}t^{\underline{k}}$ for $\alpha \in \Delta^{\times}$.   
\end{enumerate}
\ermk

\subsection{Twisted Full Toroidal Lie Algebras}\label{Twisted Toroidal}
Consider the Lie algebra of derivations $\mathcal{D}(\underline{m})$ on $A(\underline{m})$, which we can define by setting 
\begin{align}\label{Witt}
\mathcal{D}(\underline{m}) = \text{span} \{t^{\underline{r}}d_i : \underline{r} \in \Gamma, \ 1 \leqslant i \leqslant n \}.    
\end{align}
Now similar to the untwisted case, we can likewise take any linear span $\phi$ of the 2-cocyles $\phi_1$ and $\phi_2$, now defined on $\mathcal{D}(\underline{m})$, with values in $\mathcal{Z}(\underline{m})$. This subsequently gives us the corresponding twisted full toroidal Lie algebra
\begin{align*}
\widetilde{\tau}= LT \oplus \mathcal{Z}(\underline{m}) \oplus \mathcal{D}(\underline{m})
\end{align*}
by prescribing the following commutator relations along with the bracket operations (1) and (2) mentioned in Subsection \ref{Extension}.
\begin{enumerate}
	\item $[t^{\underline{r}}d_i,t^{\underline{s}}K_j]=s_it^{\underline{r} + \underline{s}}K_j + \delta_{ij} \sum_{p=1}^{n}r_p t^{\underline{r} + \underline{s}}K_p$,
	\item $[t^{\underline{r}}d_i,t^{\underline{s}}d_j]=s_it^{\underline{r} + \underline{s}}d_j-r_j t^{\underline{r} + \underline{s}}d_i + \phi(t^{\underline{r}}d_i,t^{\underline{s}}d_j)$,
	\item $[t^{\underline{r}}d_i,x \otimes t^{\underline{k}}]=k_ix \otimes t^{\underline{r} + \underline{k}} \ \forall \ x \in \mathfrak{g}_{\overline{k}},\ \underline{k} \in \mathbb{Z}^n, \ \underline{r}, \ \underline{s} \in \Gamma, \ 1 \leqslant i,j \leqslant n$.
\end{enumerate}
The centre of $\widetilde{\tau}$ is clearly spanned by the elements $K_1, \ldots, K_n$.
\brmk\label{Subalgebra} 
Observe that the universal central extension of $\mathfrak{g}_{\bar{0}} \otimes A(\underline{m})$ is simply $\mathfrak{g}_{\bar{0}} \otimes A(\underline{m}) \oplus \mathcal{Z}(\underline{m})$ (see (\ref{Universal})). It is now immediate that the (untwisted) full toroidal Lie algebra 
\begin{align*}
\tau_{A(\underline{m})} = \mathfrak{g}_{\bar{0}} \otimes A(\underline{m}) \oplus \mathcal{Z}(\underline{m}) \oplus \mathcal{D}(\underline{m})
\end{align*}
is a subalgebra of $\widetilde{\tau}$. Moreover if we take $D$=span$\{d_1, \ldots, d_n \}$, then $\widetilde{\tau}$ also contains the subalgebra, namely the $graded$ $Lie$ $torus$ given by 
\begin{align*}
\widetilde{LT}= LT \oplus \mathcal{Z}(\underline{m}) \oplus D.   
\end{align*}
\ermk

\subsection{Roots and Coroots}\label{Roots}
Consider the abelian subalgebra of $\widetilde{\tau}$ given by
\begin{align*}
\widetilde{\mathfrak{h}} = \mathfrak{h}_{\overline{0}} \oplus \displaystyle{\sum_{i =1}^{n}\CC K_i \oplus \sum_{i=1}^{n}} \CC d_i.    
\end{align*}
In order to describe the roots of $\widetilde{\tau}$, let us first define $\delta_i \in \widetilde{\mathfrak{h}}^*$ by setting
\begin{align*}
\delta_i (\mathfrak{h}_{\overline{0}}) =0,\,\,\,\  \delta_i (K_j) = 0\,\,\,\ \text{and}\,\,\,\ \delta_i (d_j) = \delta_{i j}.
\end{align*}
Put $\delta_{\underline{\beta}} = \displaystyle{\sum_{i =1}^{n}{\beta_i \delta_i}}$ for $\underline{\beta} = (\beta_1, \ldots, \beta_n) \in \mathbb{C}^n$. For $\underline{k} \in \mathbb{Z}^n$, we shall refer to the vector $\delta_{\underline{k} + \underline{\gamma}}$ as the translate of $\delta_{\underline{k}}$ by the vector $\underline{\gamma} = (\gamma_1, \ldots, \gamma_n) \in \mathbb{C}^n$.  Set $\widetilde{\Delta} = \{ \alpha + \delta_{\underline{k}}\ | \ \alpha \in \Delta_{0, \mathrm{en}}^{\times}, \ \underline{k} \in \ZZ^{n} \} \cup \{\delta_{\underline{k}} \ | \ \underline{k} \in \mathbb{Z}^n \setminus \{\underline{0}\} \}$. Then we have the root space decomposition of $\widetilde{\tau}$ with respect to $\widetilde{\mathfrak{h}}$ given by
\begin{equation}\label{Decomposition}
\widetilde{\tau} = \displaystyle{\bigoplus_{\gamma \in \widetilde{\Delta} \cup \{0\}}} {\widetilde{\tau}_{\gamma}}
\end{equation}
where 
\[
 \ {\widetilde{\tau}_{\alpha + \delta_{\underline{k}}}}= 
\begin{cases}
    \gg_{\bar{k}}(\alpha) \otimes \mathbb{C}t^{\underline{k}}, \,\,\,\,\,\,\,\,\,\,\,\,\,\,\,\,\,\,\,\,\,\,\,\,\,\,\,\,\,\,\,\,\,\,\,\,\,\,\,\,\,\,\,\,\,\,\,\,\,\,\,\,\,\,\,\,\,\,\,\,\,\,\,\,\,\,\,\,\,\,\,\,\,\ \text{if} \,\,\,\,\alpha \neq 0 \\
     \big (\mathfrak{h}_{\overline{0}} \otimes \mathbb{C}t^{\underline{k}}\big ) \oplus  \displaystyle{\big ( \sum_{i=1}^{n}}\mathbb{C}{t^{\underline{k}}}K_i \big ) \oplus  \displaystyle{\big ( \sum_{i=1}^{n}}\mathbb{C}{t^{\underline{k}}}d_i \big ) \ , \ \text{if}\,\,\,\,\alpha=0,\ \underline{k} \in \Gamma\\
     \gg_{\bar{k}}(0) \otimes \mathbb{C}t^{\underline{k}},\,\,\,\,\,\,\,\,\,\,\,\,\,\,\,\,\,\,\,\,\,\,\,\,\,\,\,\,\,\,\,\,\,\,\,\,\,\,\,\,\,\,\,\,\,\,\,\,\,\,\,\,\,\,\,\,\,\,\,\,\,\,\,\,\,\,\,\,\,\,\,\,\,\,\,\ \text{if} \,\,\,\,\alpha = 0,\ \underline{k} \notin \Gamma.    
\end{cases}
\]
This shows that the roots of $\widetilde{\tau}$ are given by $\widetilde{\Delta}$. Furthermore, let us define
\begin{align*}
\widetilde{\Delta}_{+}= \{ \alpha + \delta_{\underline{k}}\ | \ \alpha \in \Delta_{0, \mathrm{en}}^{\times,+}\ ,\ \underline{k} \in \ZZ^{n} \},\ \ \widetilde{\Delta}_{-}= \{ \alpha + \delta_{\underline{k}}\ | \ \alpha \in \Delta_{0, \mathrm{en}}^{\times,-}\ ,\ \underline{k} \in \ZZ^{n} \}
\end{align*}
as the set of all positive and negative roots of $\widetilde{\tau}$ respectively.
A root $\gamma = \alpha + \delta_{\underline{k}}$ is said to be real if $\alpha \neq 0$, else  we call it a null root. Let $\widetilde{\Delta}^{\mathrm{re}}$ be the set of all real roots of $\widetilde{\tau}$. For each $\gamma = \alpha + \delta_{\underline{k}} \in \widetilde{\Delta}^{\mathrm{re}}$, define the corresponding co-root $\gamma^{\vee} := \alpha^{\vee} + \frac{2}{(\alpha|\alpha)} \displaystyle{\sum_{i =1}^{n}{k_i K_i}}$ where $\alpha^{\vee} \in \mathfrak{h}_{\overline{0}}$ is the co-root of $\alpha \in \Delta_{0, \mathrm{en}}^{\times}$.       

\subsection{The Weyl Group}
For each $\gamma \in \widetilde{\Delta}^{\mathrm{re}}$, define
the reflection operator $r_{\gamma}$ on $\widetilde{\mathfrak{h}}^*$ by setting
\begin{align*}
r_{\gamma}(\lambda) = \lambda - \lambda(\gamma^{\vee}) \gamma\ \forall \ \lambda \in \widetilde{\mathfrak{h}}^*.
\end{align*}
Then the Weyl group of $\widetilde{\tau}$, which we shall denote by $W$ is the group
generated by all such reflections $r_{\gamma}$ with $\gamma \in \widetilde{\Delta}^{\mathrm{re}}$. 
\section{Integrable Modules}\label{Representations}
In this section, we introduce the notions of integrable representations and level zero highest weight modules. We conclude this section by ultimately showing that every irreducible integrable module over $\widetilde{\tau}$ of level zero is in fact a highest weight module with respect to a suitable triangular decomposition.
\bdfn
A $\widetilde{\tau}$-module $V$ is called integrable if 
\begin{enumerate}
\item $\displaystyle{V = \bigoplus_{\lambda \in \widetilde{\mathfrak{h}}^*} {V_{\lambda}}}$, where $V_{\lambda} = \{v \in V |\,\, h.v = \lambda(h)v \,\, \forall \,\, h \in \widetilde{\mathfrak{h}}\}$.
\item For each $x \in \mathfrak{g}_{\overline{k}}(\alpha)\otimes \mathbb{C}t^{\underline{k}}$ $(\alpha \neq 0)$ and every $v \in V$, there exists some $m=m(\alpha,\underline{k},v) \in \mathbb{N}$ such that $x^m.v=0$.      
\end{enumerate} 
\edfn

For an integrable module $V$ over $\widetilde{\tau}$, we shall denote the set of all weights of $V$ by $P(V) = \{\mu \in \widetilde{\mathfrak{h}}^*|\ V_{\mu} \neq (0)\}$. For any $\mu \in P(V)$, $V_{\mu}$ is the called the weight space of $V$ having weight $\mu$ and the elements of $V_{\mu}$ are referred to as the weight vectors of $V$ with weight $\mu$. In this paper, our goal is to classify all those irreducible integrable modules over $\widetilde{\tau}$ having finite-dimensional weight spaces with respect to $\widetilde{\mathfrak{h}}$ where the central elements $K_1, \ldots, K_n$ act trivially. These modules are commonly known as $level$ $zero$ modules.      

\blmma\label{Integrable}
Let $V$ be an integrable module for $\widetilde{\tau}$ with finite-dimensional weight spaces. Then
\begin{enumerate}
	\item $P(V)$ is invariant under the action of the Weyl group $W$.
	\item $\mathrm{dim}(V_{\lambda}) = \mathrm{dim}(V_{w \lambda})$ $\forall\ \lambda \in P(V)$ and $w \in W$.
	\item If $\lambda \in P(V)$  and $\gamma \in \widetilde{\Delta}^{\mathrm{re}}$, then $\lambda (\gamma^{\vee}) \in \ZZ$.
	\item If $\lambda \in P(V)$  and $\gamma \in \widetilde{\Delta}^{\mathrm{re}}$ with $\lambda(\gamma^{\vee}) >0$, then $\lambda - \gamma \in P(V)$.
\end{enumerate}
\elmma
\begin{proof}
The proof is similar to \cite[Lemma 2.3]{E1}. Note that the irreducibility of $V$ is not required.     
\end{proof}
The root space decomposition of $\widetilde{\tau}$ in \eqref{Decomposition} automatically induces a natural triangular decomposition of $\widetilde{\tau}$ given by
 \begin{align}\label{LT}
 \widetilde{\tau} = \widetilde{\tau}_{-} \oplus \ \widetilde{\tau}_{0} \ \oplus \widetilde{\tau}_{+}
 \end{align}
 where
\begin{align*}
& \widetilde{\tau}_{-} =  \displaystyle{ \bigoplus_{\alpha \in \Delta_{0, \mathrm{en}}^{\times,-}, \,\,\underline{k} \in \ZZ^{n}}}\bigg ({\gg_{\bar{k}}(\alpha) \otimes \mathbb{C}t^{\underline{k}}}\bigg ) ,\\
& \widetilde{\tau}_{0} = \bigg (\displaystyle{ \sum_{\underline{k} \in \ZZ^{n}}} {\gg_{\bar{k}}(0) \otimes \mathbb{C}t^{\underline{k}}} \bigg )\ \bigoplus\ \bigg (\displaystyle{ \sum_{1 \leq i \leq n, \,\,\underline{k} \in \Gamma}}\mathbb{C}{t^{\underline{k}}}K_i \bigg ) \ \bigoplus  \mathcal{D}(\underline{m}),\\
& \widetilde{\tau}_{+} =  \displaystyle{ \bigoplus_{\alpha \in \Delta_{0, \mathrm{en}}^{\times,+}, \,\,\underline{k} \in \ZZ^{n}}}\bigg ({\gg_{\bar{k}}(\alpha) \otimes \mathbb{C}t^{\underline{k}}}\bigg ).
\end{align*}

\bdfn \label{D3.3}
$V$ is said to be a level zero highest weight module for $\widetilde{\tau}$ if there exists a non-zero weight vector $v$ in $V$ such that
\begin{enumerate}
\item $V=U(\widetilde{\tau})v$.
\item $\widetilde{\tau}_{+}.v=0$.
\item $U(\widetilde{\tau}_{0})v$ is an irreducible module over $\widetilde{\tau}_{0}$.
\item The elements $K_1, \ldots, K_n$ act trivially on $V$.
\end{enumerate}
\edfn
For the rest of this section, let us fix an irreducible integrable module $V$ over $\widetilde{\tau}$ of level zero having finite-dimensional weight spaces with respect to the Cartan subalgebra $\widetilde{\mathfrak{h}}$. Further set 
\begin{align*}
V_{+}=\{v \in V \ | \ \widetilde{\tau}_{+}.v=(0)\}.
\end{align*}

\bppsn\label{Highest}
The highest weight space $V_{+}$ is non-zero.
\eppsn

\begin{proof} By Remark \ref{Subalgebra}, $V$ is clearly an integrable module over $\widetilde{LT}$ having finite-dimensional weight spaces with respect to $\widetilde{\mathfrak{h}}$. The desired result is now a direct consequence of \cite[Proposition 3.7]{S}. Also see \cite[Lemma 2.6]{E1} and \cite[Theorem 2.4(ii)]{C}. Finally observe that the irreducibility condition prescribed in \cite[Proposition 3.7]{S} is actually redundant.                
\end{proof}

\blmma\ \label{L3.5}
\begin{enumerate}
\item $V_{+}$ is an irreducible module over $\widetilde{\tau}_{0}$.
\item The weights of $V_{+}$ are the same up to a translate of a subset of the null roots, i.e. there exists a unique $\overline{\lambda} \in {{\mathfrak{h}}^*_{\overline{0}}}$ and some $\underline{\beta} \in \mathbb{C}^n$ (not necessarily unique) such that $P(V_{+}) \subseteq \{ \overline{\lambda} + \delta_{\underline{r} + \underline{\beta}} \ | \ \underline{r}  \in \mathbb{Z}^n\}$.
\end{enumerate}
\elmma

\begin{proof}
(1) The required result easily follows from the irreducibility of $V$ and by an application of the PBW theorem.\\
(2) The proof proceeds verbatim as in the case of \cite[Lemma 4.5(1)]{S}. Just observe that ${\mathfrak{h}}_{\overline{0}}$ commutes with $\mathcal{D}(\underline{m})$.
\end{proof}
\bcrlre\label{C3.6}
$V$ is a level zero highest weight module for $\widetilde{\tau}$.
\ecrlre
\begin{proof}  This is a trivial consequence of Proposition \ref{Highest} and Lemma \ref{L3.5}.
\end{proof}

\brmk\label{R3.7}
If $\overline{\lambda}=0$ in Lemma \ref{L3.5}, then using the integrability of $V$, it is easy to conclude that the only possible weights of $V$ are $\delta_{\underline{\gamma}}$ $(\underline{\gamma} \in \mathbb{C}^n)$. Then from Remark \ref{Generate}, it directly follows that $\overline{LT}$ acts trivially on $V$. This ultimately reduces our problem to classifying all the irreducible modules for the algebra $\mathcal{D}(\underline{m})$ having finite-dimensional weight spaces with respect to the space of degree derivations $D$=span$\{d_1, \ldots, d_n \}$. The classification of these irreducible modules can be found in \cite{BF}. In the present paper, we shall primarily focus on those irreducible modules where $\overline{\lambda} \neq 0$.
\ermk
\section{Action Of The Central Operators}\label{Central}
Throughout this section, unless otherwise explicitly stated, $V$ will always stand for a level zero irreducible (but not necessarily integrable) highest weight module for $\widetilde{\tau}$ (see Definition \ref{D3.3}) having finite-dimensional weight spaces with respect to $\widetilde{\mathfrak{h}}$. We shall also assume that the Cartan subalgebra ${{\mathfrak{h}}_{\overline{0}}}$ of ${{\mathfrak{g}}_{\overline{0}}}$ does not act trivially on the highest weight space $V_{+}$. Under these assumptions, our main aim in this section is to prove the following result.
\bthm\label{T4.1}
$\mathcal{Z}(\underline{m})$ acts trivially on $V$.
\ethm
\noindent We need some preparation to prove this theorem. First observe that
by the initial assumptions on our module $V$ and Lemma \ref{L3.5}, there exists  $\Lambda \in P(V)$ with $\Lambda|_{\mathfrak{h}_{\overline{0}}} \neq 0$ satisfying
\begin{align}\label{Recall}
h.w = \Lambda(h)w \ \forall \ h \in {\mathfrak{h}}_{\overline{0}},\ w \in V_{+}.
\end{align}
Pick some $h_0 \in \mathfrak{h}_{\overline{0}}$ such that $\Lambda(h_0) \neq 0$. It is now clear that
\begin{align*}
V_{+} = \bigoplus_{\underline{m} \in \mathbb{Z}^n} V_{+}(\underline{m})
\end{align*}
with $V_{+}(\underline{m}) = \{v \in V_{+}\ | \ d_i.v = (\Lambda(d_i)+m_i)v, \ 1 \leqslant i \leqslant n \}$ which shows that $V_{+}$ is a $\mathbb{Z}^n$-graded $\widetilde{\tau}_0$-module with finite-dimensional graded components.  

\blmma\label{L4.2}
Let $h_0 \otimes t^{\underline{k}} \in \widetilde{\tau}_0$ where $\underline{k} \in \Gamma \setminus \{\underline{0}\}$. Suppose that there exists a non-zero element $w \in V_{+}$ such that $(h_0 \otimes t^{\underline{k}}).w=0$. Then $h_0 \otimes t^{\underline{k}}$ is locally nilpotent on $V_{+}$.  
\elmma

\begin{proof}
It is trivial to check that for any $\underline{r} \in \Gamma$ and $i = 1, \ldots, n$, we have 
\begin{align*}
{(h_0 \otimes t^{\underline{k}})}^3 ((t^{\underline{r}}d_i).w) = 0.
\end{align*}
Similarly we can show that if $x_{\overline{r}} \otimes t^{\underline{r}} \in \big(\mathfrak{g}_{\overline{r}}(0) \otimes \mathbb{C}t^{\underline{r}} \big )$ where $\underline{r} \in \mathbb{Z}^n$, then
\begin{align*}
{(h_0 \otimes t^{\underline{k}})}^2 ((x_{\overline{r}} \otimes t^{\underline{r}}).w) = 0.
\end{align*}
Consequently by induction on $p$, it can be further deduced that
\begin{align*}
{(h_0 \otimes t^{\underline{k}})}^{2p+1} \big ((t^{\underline{r_1}}d_{i_1} \ldots t^{\underline{r_p}}d_{i_p}).w \big ) = 0 \ ,\\
{(h_0 \otimes t^{\underline{k}})}^{p+1} \bigg (\big ((x_{\overline{s_1}} \otimes t^{\underline{s_1}}) \ldots (x_{\overline{s_p}} \otimes t^{\underline{s_p}}) \big ) .w \bigg ) = 0
\end{align*}
for all $p \in \mathbb{N}$, $\underline{r_1}, \ldots, \underline{r_p} \in \Gamma$ and $\underline{s_1}, \ldots, \underline{s_p} \in \mathbb{Z}^n$ with $1 \leqslant i_1, \ldots, i_p \leqslant n$.\\
Finally since $\mathfrak{h}_{\overline{0}} \otimes A(\underline{m})$ commutes with $\mathcal{Z}(\underline{m})$, the lemma now directly follows from the irreducibility of $V_{+}$ over $\widetilde{\tau}_0$.   
\end{proof}

Next let us record some results which can be proved using the above Lemma \ref{L4.2} and Remark \ref{Subalgebra} in a more or less similar manner as in \cite{CD} and \cite{EJ} by considering $h_0 \otimes t^{\underline{k}}$ instead of $t^{\underline{m}}k_0$ and replacing $c_0$ by $\Lambda(h_0)$. Our proof of Theorem \ref{T4.1} is independent of the final assertion (Lemma \ref{Finite}(4)) and this assertion can be deduced by proceeding exactly as in \cite{CD,EJ} after applying our Theorem \ref{T4.1}.      
\blmma\label{Finite}\
\begin{enumerate}
\item dim $V_{+}(\underline{k}) = $ dim $V_{+}(\underline{k} + \underline{r})=N_{\underline{k}} \in \mathbb{N} \ \forall \ \underline{k} \in \mathbb{Z}^n, \ \underline{r} \in \Gamma$. In particular, we have dim $V_{+}(\underline{r}) = $ dim $V_{+}(\underline{s})=N$(say) $\forall \ \underline{r} , \underline{s} \in \Gamma$.       
\item Let $w$ be a non-zero element in $V_{+}$ such that $(t^{\underline{r}}K_i).w=0$ for some $\underline{r} \in \Gamma \setminus \{\underline{0}\}$ and $1 \leqslant i \leqslant n$. Then $t^{\underline{r}}K_i$ acts locally nilpotently on $V_{+}$.
\item For each $\underline{r} \in \Gamma, \ (t^{\underline{r}}K_i)^N.V_{+}=(0) \ \forall \ i=1, \ldots, n$.
\item $\big ( (h_0 \otimes t^{\underline{r}}) (h_0 \otimes t^{\underline{s}}) \big ).v= \Lambda(h_0)(h_0 \otimes t^{\underline{r}+ \underline{s}}).v \ \forall \ \underline{r}, \underline{s} \in \Gamma, \ v \in V_{+}$.   
\end{enumerate}
\elmma

\blmma\label{L4.4}  
Let $n \geqslant 2$ and suppose that for each $1 \leqslant i \leqslant n$, there exists a fixed $N \in \mathbb{N}$ such that $(t^{\underline{r_1}}K_i \ldots t^{\underline{r_N}} K_i).V_{+}=(0) \ \forall \ \underline{r_1}, \ldots, \underline{r_N} \in \Gamma \setminus \{\underline{0} \}$. Then there exists a non-zero vector $v_0$ in the highest weight space $V_{+}$ such that
$\mathcal{Z}(\underline{m}).v_0=0$.  
\elmma
\begin{proof}
The proof of this lemma is algorithmic, i.e. we shall provide a precise procedure that will enable us to acquire a common eigenvector. First note that if $N=1$, then we are done. So assume that $N>1$. Now by our hypothesis, we can find a non-zero $v \in V_{+}$ satisfying $t^{\underline{r}}K_1.v=0 \ \forall \ \underline{r} \in \Gamma$. If $t^{\underline{s}}K_2.v \neq 0$ for some  $\underline{s} \in \Gamma$, then repeatedly act central operators of the form $t^{\underline{p}}K_2$ to  this vector until we obtain a non-zero vector of the form $w= (t^{\underline{r_1}}K_2 \ldots t^{\underline{r_m}} K_2).v$ for some $m \in \mathbb{N}$ such that $t^{\underline{p}}K_2.w=0 \ \forall \ \underline{p} \in \Gamma$. Our hypothesis guarantees that such a vector indeed exists and in fact $m \leqslant N-1$. This clearly gives $t^{\underline{r}}K_1.w=0=t^{\underline{r}}K_2.w=0 \ \forall \ \underline{r} \in \Gamma$. Continuing this process for the remaining operators $t^{\underline{r}}K_i$, we finally get our common eigenvector.     
\end{proof}

\noindent \textbf{Proof of Theorem \ref{T4.1}}. If $n=1$, then the assertion clearly holds good. Let us first fix any $1 \leqslant i \leqslant n$. Since we are assuming $n \geqslant 2$, we can now choose $1 \leqslant j \leqslant n$ such that $i \neq j$. Then using Lemma \ref{Finite}, we have 
\begin{align*}
0=t^{\underline{r_1}}d_j \big ((t^{\underline{s}}K_i)^N.v \big )=Ns_j \bigg (t^{\underline{r_1}+ \underline{s}}K_i \big ((t^{\underline{s}}K_i)^{N-1}.v \big ) \bigg ) \ \forall \ \underline{r_1}, \underline{s} \in \Gamma, \ v \in V_{+}
\end{align*}
from which it clearly follows that
\begin{align*}
\big (t^{\underline{r_1}+ \underline{s}}K_i (t^{\underline{s}}K_i)^{N-1} \big ).V_{+}=(0) \ \forall \ \underline{r_1} \in \Gamma \ \text{and} \ \underline{s} \in \Gamma \ \text{satisfying} \ s_j \neq 0. 
\end{align*}
This consequently implies that 
\begin{align*}
t^{\underline{r_2}}d_j \bigg (s_jt^{\underline{r_1}+ \underline{s}}K_i \big ((t^{\underline{s}}K_i)^{N-1}.v \big ) \bigg )=0 \ \forall \ \underline{r_1}, \underline{r_2}, \underline{s} \in \Gamma, \ v \in V_{+}
\end{align*}
from which we directly get
\begin{align*}
(N-1)s_{j}^{2} \bigg (t^{\underline{r_1}+ \underline{s}}K_it^{\underline{r_2}+ \underline{s}}K_i  \big ((t^{\underline{s}}K_i)^{N-2}.v \big ) \bigg )=0
\end{align*}
for all $\underline{r_1},\ \underline{r_2} \in \Gamma$ with $\underline{s} \in \Gamma$ satisfying $s_j \neq 0$ and $v \in V_{+}$. Then
\begin{align*}
\big (t^{\underline{r_1}+ \underline{s}}K_it^{\underline{r_2}+ \underline{s}}K_i (t^{\underline{s}}K_i)^{N-2} \big ). V_{+}=(0) \ \forall \ \underline{r_1},\ \underline{r_2} \in \Gamma \ \text{and} \ \underline{s} \in \Gamma \ \text{with} \ s_j \neq 0.
\end{align*}
By repeating the above argument another $(N-2)$ times, we thus obtain
\begin{align*}
(t^{\underline{r_1}+ \underline{s}}K_i \ldots t^{\underline{r_N}+ \underline{s}}K_i). V_{+}=(0)
\end{align*}
for all $\underline{r_1}, \ldots, \underline{r_N} \in \Gamma$ and $\underline{s} \in \Gamma$ satisfying $s_j \neq 0$. Let us consider arbitrary $\underline{s_1}, \ldots, \underline{s_N} \in \Gamma \setminus \{\underline{0} \}$. Choose any $\underline{s} \in \Gamma$ such that $s_j \neq 0$ and thereby set $\underline{r_k}= \underline{s_k}- \underline{s} \ \forall \ k=1, \ldots, N$. This finally gives us
\begin{align*}
(t^{\underline{s_1}}K_i \ldots t^{\underline{s_N}}K_i). V_{+}=(0) \ \forall \ i=1, \ldots, n.
\end{align*}
Subsequently by Lemma \ref{L4.4}, there exists a non-zero vector $v_0 \in V_{+}$ such that $(t^{\underline{r}}K_i).v_0=0 \ \forall \ \underline{r} \in \Gamma, \ 1 \leqslant i \leqslant n$. Now it is trivial to check that
\begin{align*}
W= \{v \in V \ | \ \mathcal{Z}(\underline{m}).v=0 \} 
\end{align*}
is a $\widetilde{\tau}$-submodule of $V$. Therefore we are done by the irreducibility of $V$.
\brmk
The above theorem thereby allows us to fully restrict our attention to level zero irreducible modules over $\widehat{\tau}$ having finite-dimensional weight spaces with respect to $\widehat{\mathfrak{h}}$  where $\widehat{\tau}= LT \oplus \mathcal{D}(\underline{m})$ and $\widehat{\mathfrak{h}} = \mathfrak{h}_{\overline{0}} \oplus D$.   
\ermk
\noindent The next result follows from Corollary \ref{C3.6}, Remark \ref{R3.7} and Theorem \ref{T4.1}.       
\bcrlre\label{C4.9}
Let $V$ be a level zero irreducible integrable module for $\widetilde{\tau}$ having finite-dimensional weight spaces with respect to $\widetilde{\mathfrak{h}}$. Then $\mathcal{Z}(\underline{m})$ acts trivially on $V$.  
\ecrlre

\brmk\
\begin{enumerate}
\item The proof of Corollary \ref{C4.9} for the $untwisted$ full toroidal Lie algebras provided in \cite[Theorem 4.1]{EJ} is not complete. This is because unlike in $\mathcal{Z}(\underline{m})$, the operators $h_0 \otimes t^{\underline{r}}$ a priori $do$ $not$ commute with each other. Consequently we cannot employ the techniques suggested by the authors in \cite{EJ} to get a commutative Lie algebra (denoted by $\mathscr{C}$ in \cite[Theorem 3.1]{EJ}) in the level zero case. Nevertheless, their theorem follows instantaneously from our Theorem \ref{T4.1}.    
\item The above corollary can be also derived independently based on the arguments sketched in \cite[Theorem 4.9]{S}. But we chose a different approach in this case as it helps us to conclude that the assertion actually holds good for a much bigger class of level zero highest weight modules that may not be integrable.
\end{enumerate}
\ermk

\section{Finite-Dimensional Graded-Irreducible Representations and Evaluation Modules}\label{Graded-Irreducible}
This section is mainly devoted to the recollection of a few notions and results from \cite{EK} and \cite{L} which will be useful for our classification problem. At the end of this section, we shall prove
an important result that will be utilized in the next section.
\subsection{Modules twisted by an automorphism}\cite[Subsection 3.1]{EK}\label{SS5.1}\\ 
 Recall that $\mathfrak{g}$ is a finite-dimensional simple Lie algebra also equipped with a $G$-grading. Putting $\widehat{G}= \{f:G \longrightarrow \mathbb{C}^{\times} \ | \ f \text{ is a group homomorphism} \}$ which is a finite group, it is immediate that any $\chi \in \widehat{G}$ gives rise to an automorphism of $\mathfrak{g}$ defined by
\begin{align*}
\alpha_{\chi}(x)= \chi(\overline{k})x \ \forall \ x \in \mathfrak{g}_{\overline{k}}, \ \underline{k} \in \mathbb{Z}^n.  
\end{align*}
For any $\mathfrak{g}$-module $V$ and $\phi \in$ Aut($\mathfrak{g}$), we also have the twisted $\mathfrak{g}$-module $V^{\phi}$, which is nothing but the same vector space $V$, but now twisted by a $\mathfrak{g}$-action given by $x.v= \phi(x)v$ for all $x \in \mathfrak{g}$ and $v \in V$. Thus we find that the group Aut($\mathfrak{g}$) and subsequently $\widehat{G}$ acts (on the right) on the isomorphism classes of $\mathfrak{g}$-modules which thereby induces a natural action of $\widehat{G}$ on $P_{\mathfrak{g}}^{+}$. Consequently the action of $\widehat{G}$ on the isomorphism classes of finite-dimensional irreducible $\mathfrak{g}$-modules can be represented as an action of the Dynkin automorphisms of $\mathfrak{g}$ on the set of dominant integral weights of $\mathfrak{g}$ by simply permuting the vertices of the Dynkin diagram and the corresponding fundamental weights. Let us denote this action by `$\circ$'. The finite-dimensional irreducible highest weight module for $\mathfrak{g}$ with highest weight $\lambda$ will be denoted by $V(\lambda)$.  

\subsection{Finite-dimensional graded-irreducible modules}\label{SS5.2}
Let us consider $\lambda \in P_{\mathfrak{g}}^{+}$ with $graded \ Schur \ index$ S$(\lambda)\in \mathbb{N}$ (see \cite[Section 3]{EK} for the precise definition of S($\lambda$)). If $\widehat{G} \lambda= \{ \lambda_1, \ldots, \lambda_l \}$, then it can be shown that $\bigoplus_{i=1}^{l} V(\lambda_i)^{\oplus S(\lambda)}$ is a finite-dimensional $G$-graded-irreducible module for $\mathfrak{g}$ (see \cite[Subsection 3.3]{EK}). Conversely, up to an isomorphism of $(ungraded)$ modules over $\mathfrak{g}$, every finite-dimensional $G$-graded-irreducible module for $\mathfrak{g}$ is of the form $\bigoplus_{i=1}^{l} V(\lambda_i)^{\oplus S(\lambda)}$ where $\widehat{G} \lambda= \{ \lambda_1, \ldots, \lambda_l \}$ for some $\lambda \in P_{\mathfrak{g}}^{+}$ (see \cite[Theorem 8]{EK}). Moreover it is also known that $\widehat{G} \lambda_i= \widehat{G} \lambda_j$ and $S(\lambda_i)=S(\lambda_j) \ \forall \ 1 \leqslant i,j \leqslant l$ (see \cite[Subsection 3.4]{EK}). We finally remark that if $G$ is trivial, then we always have S($\lambda$)=1.   
\subsection{Finite-dimensional evaluation modules} For any $\underline{a} \in ({\mathbb{C}^{\times}})^{n}$ and $\underline{k} \in \mathbb{N}^{n}$, write $\underline{a}^{\underline{k}}= \prod_{i=1}^{n}a_i^{k_i}$ and  
$\underline{a}(\underline{k})=(a_1^{k_1}, \ldots, a_n^{k_n})$. Consider $\underline{m} \in \mathbb{N}^{n}$ as mentioned in Subsection \ref{Twisted}. We say that $\underline{a} \in \mathbb{C}^{n}$ is a root of unity for $\underline{m}$ if we have
$\underline{a}(\underline{m})=(1, \ldots, 1)= \underline{1}$. Let $\mathcal{U}(\underline{m})$ represent the set of all roots of unity for $\underline{m}$ in $\mathbb{C}^n$. We can now put an $LT$-module structure on the irreducible $\mathfrak{g}$-module $V(\lambda)$ by setting 
\begin{align*}
(x \otimes t^{\underline{k}}).v= \underline{a}^{\underline{k}} (x.v) \ \forall \ x \in \mathfrak{g}_{\overline{k}} \ , \ \underline{k} \in \mathbb{Z}^n, \ v \in V(\lambda). 
\end{align*}
We shall denote this $LT$-module by $ev_{\underline{a}}V(\lambda)$.
These modules are commonly known as $evaluation \ modules$ in the literature which can be easily shown to be irreducible (see \cite[Section 4]{L} for more details). The finite-dimensional irreducible $LT$-modules have been completely classified in \cite{L} and \cite{NSS}.   

\bppsn\cite[Theorem 5.4]{L}\label{P5.1}
The finite-dimensional $LT$-modules $ev_{\underline{a}}V(\lambda)$
and $ev_{\underline{b}}V(\mu)$ are isomorphic if and only if we have $\underline{a}(\underline{m})=\underline{b}(\underline{m})$ and $\mu= \lambda \circ \gamma$ where $\gamma$ is the outer part of $\omega \in$ Aut($\mathfrak{g}$) that is defined by $\omega(x)= (\underline{a}^{\underline{k}}/\underline{b}^{\underline{k}})x$ for all $\underline{k} \in \mathbb{Z}^n$ and $x \in \mathfrak{g}_{\overline{k}}$ .    
\eppsn

\bcrlre\label{C5.2}
Let $\lambda, \mu \in P_{\mathfrak{g}}^{+}$ with $\mu \in \widehat{G}\lambda$. Then for any given $\underline{a} \in ({\mathbb{C}^{\times}})^{n}$, there exists some $\underline{b} \in ({\mathbb{C}^{\times}})^{n}$ such that $ev_{\underline{b}}V(\lambda) \cong ev_{\underline{a}}V(\mu)$ as $LT$-modules. 
\ecrlre
\begin{proof}
We know that $G \cong G_1 \times \ldots \times G_n$ where $G_i= \langle{g_i}\rangle$ (say) $\forall \ 1 \leqslant i \leqslant n$. By 
hypothesis, there exists a Dynkin automorphism $\gamma$ satisfying $\mu= \lambda \circ \gamma$. Moreover $\gamma$ is the outer part of some $\omega \in$ Aut($\mathfrak{g}$) given by
\begin{align*}
\omega(x)= \chi(\overline{k})x \ \forall \ \underline{k}=(k_1, \ldots, k_n) \in \mathbb{Z}^n, \ x \in \mathfrak{g}_{\overline{k}}
\end{align*}	
where $\chi \in \widehat{G}$. Now since the group $\widehat{G}$ is canonically isomorphic to the group $\widehat{G_1} \times \ldots \times \widehat{G_n}$, we can use this identification to get a unique element 
$(\chi_1, \ldots, \chi_n) \in \widehat{G_1} \times \ldots \times \widehat{G_n}$ corresponding to our $\chi \in \widehat{G}$. Moreover as $|G_i|=m_i$, there exists some $d_i \in \mathbb{Z}$ such that $\chi_i(g_i)= \xi_{i}^{d_i}$, where $\xi_i$ is a $m_i-$th primitive root of unity. Finally setting $b_i= \xi_{i}^{d_i}a_i$, we thereby obtain
\begin{align*}
w(x)= \bigg (\prod_{i=1}^{n} \chi_i(\overline{k_i}) \bigg )x = (\underline{b}^{\underline{k}}/\underline{a}^{\underline{k}})x 
\end{align*} 	
for all $\underline{k} \in \mathbb{Z}^n$ and $x \in \mathfrak{g}_{\overline{k}}$, where $\underline{b}=(b_1, \ldots, b_n) \in 
({\mathbb{C}^{\times}})^{n}$. The assertion is now a direct consequence of Proposition \ref{P5.1}. 	
\end{proof}

\section{An Explicit Realization}\label{Realization}
In this section, we shall extract a family of level zero integrable modules for $\widetilde{\tau}$ having finite-dimensional weight spaces with respect to the Cartan subalgebra $\widetilde{\mathfrak{h}}$ and eventually show that they are all irreducible.

Let $\psi \in P_{\mathfrak{sl}_n}^{+}$ and $c \in \mathbb{C}$. It is well-known that these parameters give rise to a unique finite-dimensional irreducible $\mathfrak{gl}_n(\mathbb{C})$-module, say $V_1 = V(c, \psi)$. Again take any $\lambda \in P_{\mathfrak{g}}^{+} \setminus \{0\}$ (which we shall henceforth denote by $(P_{\mathfrak{g}}^{+})^{\times}$) with $\widehat{G} \lambda= \{ \lambda_1, \ldots, \lambda_l \}$ and graded Schur index $S(\lambda)=k$ (say). Then from our discussion in Subsection \ref{SS5.2}, it is evident that $V_{2}^{\prime}= \bigoplus_{i=1}^{l} V(\lambda_i)^{\oplus{k}}$ is a finite-dimensional $G$-graded-irreducible $\mathfrak{g}$-module. Observe that the number of irreducible summands occurring in the above decomposition of $V_{2}^{\prime}$ is equal to $kl=N$ (say). Take $V_2= \bigoplus_{\nu=1}^{N}V_{2}^{\nu}$, where the irreducible components $V_{2}^{\nu}$ of $V_2$ are all coming from the decomposition of $V_{2}^{\prime}$. For example, if $k=l=2$, there are precisely five choices for our $V_2$. Now since $\widehat{G} \lambda= \{ \lambda_1, \ldots, \lambda_l \}$ and $\widehat{G}\lambda_i= \widehat{G}\lambda_j$ for all $1 \leqslant i,j \leqslant l$, we can invoke Corollary \ref{C5.2} to pick evaluation points $\underline{a_1}, \ldots, \underline{a_N} \in \mathcal{U}(\underline{m})$ such that the twisted $\mathfrak{g}$-module (which we shall again denote by $V_2$) defined by 
\begin{align}\label{Ungraded} 
x. \big (\sum_{\nu=1}^{N} v_{2}^{\nu} \big)= \sum_{\nu=1}^{N} \underline{a_{\nu}}^{\underline{l}}(x.v_{2}^{\nu}) \ \forall \ x \in \mathfrak{g}_{\overline{l}}, \ v_{2}^{\nu} \in V_{2}^{\nu}
\end{align}
is isomorphic to $V_{2}^{\prime}$ as ungraded $\mathfrak{g}$-modules. This immediately suggests that this twisted module is $G$-graded-irreducible over $\mathfrak{g}$.\\
Let $\{E_{ij} \ | \ 1 \leqslant i,j \leqslant n \}$ be the standard basis of $\mathfrak{gl}_n(\mathbb{C})$. Now pick any $\underline{\beta}=(\beta_1, \ldots, \beta_n) \in \mathbb{C}^n$ and define an action of $\widetilde{\tau}$ on $V_1 \otimes V_2 \otimes A$ by setting
\begin{align*}
(x \otimes t^{\underline{l}}). \big (\sum_{\nu=1}^{N}v_1 \otimes v_{2}^{\nu} \otimes t^{\underline{k}} \big )= \sum_{\nu=1}^{N}v_1 \otimes \underline{a_{\nu}}^{\underline{l}}(x.v_{2}^{\nu}) \otimes t^{\underline{k} + \underline{l}} \ ;\\
t^{\underline{r}}d_i.\big (\sum_{\nu=1}^{N}v_1 \otimes v_{2}^{\nu} \otimes t^{\underline{k}} \big )=(k_i + \beta_i) \big (\sum_{\nu=1}^{N}v_1 \otimes v_{2}^{\nu} \otimes t^{\underline{k} + \underline{r}} \big )\\
+ \ \sum_{\nu=1}^{N} \sum_{j=1}^{n} \big ((r_jm_jE_{ji}).v_1 \big ) \otimes v_{2}^{\nu} \otimes t^{\underline{k} + \underline{r}} \ ;\\
t^{\underline{r}}K_i.\big (\sum_{\nu=1}^{N}v_1 \otimes v_{2}^{\nu} \otimes t^{\underline{k}} \big )=0
\end{align*}
for all $x \in \mathfrak{g}_{\overline{l}}, \ \underline{k}, \ \underline{l} \in \mathbb{Z}^n, \ \underline{r} = (r_1m_1, \ldots, r_nm_n) \in \Gamma, \ v_1 \in V_1, \ v_{2}^{\nu} \in V_{2}^{\nu}$ and $1 \leqslant i \leqslant n$. It can be verified that $V_1 \otimes V_2 \otimes A$ is a $\widetilde{\tau}$-module. Taking $V_2=  \bigoplus_{\overline{k} \in G} V_{2, \overline{k}}$ , let us further define $V^{\prime}:= \bigoplus_{\underline{k} \in \mathbb{Z}^n}V_1 \otimes V_{2, \overline{k}} \otimes \mathbb{C}t^{\underline{k}}$ which is again a submodule of $V_1 \otimes V_2 \otimes A$.   

\bppsn\label{P6.1}  
$V^{\prime}$ is a level zero irreducible integrable $\widetilde{\tau}$-module with finite-dimensional weight spaces having a non-trivial $LT$-action.   
\eppsn

 \begin{proof}
By our construction, the twisted module $V_2$ is a $G$-graded-irreducible module over $\mathfrak{g}$ which implies that $\bigoplus_{\underline{k} \in \mathbb{Z}^n} V_{2, \overline{k}} \otimes \mathbb{C}t^{\underline{k}}$ is an irreducible module over $LT \oplus D$ where $D$=span$\{d_1, \ldots, d_n \}$. The desired result can be now deduced using \cite[Proposition 2.8]{E2} which conveys that $V_1 \otimes A(\underline{m})$ is an irreducible module for $\mathcal{D}(\underline{m}) \rtimes \big (\mathfrak{h}_{\overline{0}} \otimes  A(\underline{m}) \big )$. Finally observe that the non-zero $\lambda \in P_{\mathfrak{g}}^{+}$ ensures us that $LT$ $does$ $not$ act trivially on $V^{\prime}$. 
\end{proof}
 
\brmk\label{R6.2}
The irreducible module $V^{\prime}$ over $\widetilde{\tau}$ is completely determined by the quadruplet $(\psi,c, \lambda, \underline{\beta}) \in P_{\mathfrak{sl}_n}^{+} \times \mathbb{C} \times (P_{\mathfrak{g}}^{+})^{\times} \times \mathbb{C}^n$ together with the corresponding multiplicities, say $p_1, \ldots, p_l$ of $V(\lambda_1), \ldots, V(\lambda_l)$ respectively that add up to $N$ ($p_i$ is allowed to be $0$) and some suitably chosen evaluation points $\underline{a_1}, \ldots, \underline{a_N} \in \mathcal{U}(\underline{m})$ where $N=kl$, $k=S(\lambda)$ and $l=|\widehat{G} \lambda|$. So we shall denote this module $V^{\prime}$ by $V_{\lambda, \underline{\beta}}^{c, \psi}(\underline{a_1}, \ldots, \underline{a_N};p_1, \ldots, p_l)$. Now for any arbitrary choice of non-negative integers $p_1, \ldots, p_l$ satisfying $\sum_{i=1}^{l}p_i=N$, we have $V_{\lambda, \underline{\beta}}^{c, \psi}(\underline{a_1}, \ldots, \underline{a_N};p_1, \ldots, p_l) \cong V_{\lambda, \underline{\beta}}^{c, \psi}(\underline{1}, \ldots, \underline{1};k, \ldots, k)$ as irreducible $\widetilde{\tau}$-modules where 
$\underline{a_1}, \ldots, \underline{a_N} \in \mathcal{U}(\underline{m})$ are chosen in a way (may not be unique) such that the twisted $\mathfrak{g}$-module $V_2$ (see (\ref{Ungraded})) used in our construction of $V_{\lambda, \underline{\beta}}^{c, \psi}(\underline{a_1}, \ldots, \underline{a_N};p_1, \ldots, p_l)$ becomes $G$-graded-irreducible. Therefore we can simply drop the evaluation points and denote these irreducible modules over $\widetilde{\tau}$ by just $V(\psi, c, \lambda, \underline{\beta})$. Finally note that $S(\lambda)$ and $\widehat{G} \lambda$ are automatically fixed after our initial choice of $\lambda$.  
\ermk

\section{Classification of Irreducible Modules}\label{Classification}
In this section, $V$ will always denote an irreducible integrable level zero module for $\widetilde{\tau}$ with finite-dimensional weight spaces and having a non-trivial action of $LT$. We shall finally show that these irreducible integrable modules are essentially of the form $V(\psi, c, \lambda, \underline{\beta})$ as described in the previous section.
                   
Recall that $\widehat{\tau}= LT \oplus \mathcal{D}(\underline{m})$ which like $\widetilde{\tau}$ has a triangular decomposition given by
$\widehat{\tau}= \widehat{\tau}_{-} \oplus \widehat{\tau}_{0} \oplus \widehat{\tau}_{+}$ where $\widehat{\tau}_{\pm}= \widetilde{\tau}_{\pm}$ and $\widehat{\tau}_{0}= \widetilde{\tau}_{0} \setminus \mathcal{Z}(\underline{m})$. Set
\begin{align*}  
\fma= \{x \in \mathfrak{g} \ | \ [h,x]=0 \ \forall \ h \in \mathfrak{h}_{\overline{0}} \}.
\end{align*}
Now since $\fma$ is invariant under the $\sigma_i$'s, this thereby induces a $G$-grading on $\fma$ with $\fma_{\overline{k}}= \mathfrak{g}_{\overline{k}} \cap \fma \ \forall \ \overline{k} \in G$. Denote the corresponding multiloop algebra by
\begin{align*}
L(\fma, \sigma)= \bigoplus_{\underline{k} \in \mathbb{Z}^n} \fma_{\overline{k}} \otimes \mathbb{C}t^{\underline{k}} \ .    
\end{align*}

We now introduce a particular Lie subalgebra of  $\mathcal{D}(\underline{m})$ which will play an important role not only in this section, but also in the rest of this paper. Let $D(\underline{u}, \underline{r})= \sum_{i=1}^{n}u_it^{\underline{r}}d_i$ and $I(\underline{u}, \underline{r})= D(\underline{u}, \underline{r})-D(\underline{u}, \underline{0}) \ \forall \ \underline{u} \in \mathbb{C}^n, \ \underline{r} \in \Gamma$.\\
It can be easily verified that for $\underline{u}, \underline{v} \in \mathbb{C}^n$ and $\underline{r}, \underline{s} \in \Gamma$,
\begin{align*}
[I(\underline{u}, \underline{r}), I(\underline{v}, \underline{s})]=(\underline{v},\underline{r})I(\underline{u},\underline{r})-(\underline{u},\underline{s})I(\underline{v},\underline{s})+I(\underline{w},\underline{s}+\underline{r})
\end{align*}
where $(\cdot,\cdot)$ is the standard inner product on $\mathbb{C}^n$ and $\underline{w}=(\underline{u}, \underline{s})\underline{v}- (\underline{v}, \underline{r}) \underline{u}$.\\
This immediately reveals that the subspace 
\begin{align*}
I:= \text{span} \{I(\underline{u}, \underline{r}) : \underline{u} \in \mathbb{C}^n, \ \underline{r} \in \Gamma \}
\end{align*}
forms a subalgebra of $\mathcal{D}(\underline{m})$.
Furthermore let us also set
\begin{align*}
\widetilde{L} = I \ltimes L(\fma, \sigma) \ \text{and} \ W= \text{span} \{(h_0 \otimes t^{\underline{r}}).v - v \ | \ \underline{r} \in \Gamma, \ v \in V_{+} \}. 
\end{align*} 
Now since $\Lambda(h_0) \neq 0$, we may assume that $\Lambda(h_0)=1$ (see (\ref{Recall})) and so by Lemma \ref{Finite}, it easily follows that $W$ is an $\widetilde{L}$-module. This again implies that $\widetilde{V} := V_{+}/W$ is also an $\widetilde{L}$-module. Finally observe that
\begin{align}\label{Grading} 
\widehat{\tau}_0 = L(\fma, \sigma) \oplus \mathcal{D}(\underline{m}), \   L(\fma, \sigma) = \bigg (\displaystyle{ \sum_{\underline{k} \in \ZZ^{n}}} {\gg_{\bar{k}}(0) \otimes \mathbb{C}t^{\underline{k}}} \bigg ).
\end{align}
\brmk
Note that $\widetilde{L}$ is naturally $G$-graded with $\widetilde{L}_{\overline{0}}= I \ltimes \mathfrak{h}_{\overline{0}} \otimes A(\underline{m})$. Moreover the $\mathbb{Z}^n$-grading on $V_{+}$ again gives rise to a $G$-grading on both $V_{+}$ and $W$ via the natural map given in (\ref{Natural}). This in turn shows us that the quotient $\widetilde{V}$ is also $G$-graded.  
\ermk
\blmma\
\begin{enumerate}
\item $W$ is a proper $\widetilde{L}$-submodule of $V_{+}$.
\item $\widetilde{V}$ is a finite-dimensional $\widetilde{L}$-module.
\end{enumerate}
\elmma
\begin{proof}
Let $z_i=h_0 \otimes t_i^{m_i}$ for each $i=1, \ldots, n$. Since $\Lambda(h_0)=1$, it can be easily deduced from Lemma \ref{Finite} that 
\begin{align*}
W= \text{span} \{z_i.v - v \ | \ v \in V_{+}, \ 1 \leqslant i \leqslant n \}.
\end{align*}
(1) With regard to the above discussion, the proof of the assertion follows verbatim as in \cite[Proposition 5.4(3)]{S}.\\
(2) From Lemma \ref{Finite}, it is evident that all the $z_i$'s are invertible on $V_{+}$. The desired result can now be established by proceeding similarly as in Claim 2 and Claim 3 in \cite[Theorem 4.5]{E1} and then applying our Lemma \ref{L3.5}. 
\end{proof}

Consider $\underline{\beta} \in \mathbb{C}^n$ as obtained in Lemma \ref{L3.5}. Then for any $\widetilde{L}$-module $V_1$, we can now put a $\widehat{\tau}_{0}$-module structure on $L(\underline{\beta},V_1)=V_1 \otimes A$ by setting
\begin{align}\label{Action}
x \otimes t^{\underline{k}}.(v_1 \otimes t^{\underline{s}}) & =  \big ((x \otimes t^{\underline{k}}).v_1 \big ) \otimes t^{\underline{k}+\underline{s}} \ ,\\
D(\underline{u},\underline{r}).(v_1 \otimes t^{\underline{s}}) & =   (I(\underline{u}, \underline{r}).v_1) \otimes t^{\underline{r}+ \underline{s}} + (\underline{u}, \underline{s}+ \underline{\beta})(v_1 \otimes t^{\underline{s} + \underline{r}})
\end{align}
for all $v_1 \in V_1, \ x \in \fma_{\overline{k}}$ with $\underline{k}, \ \underline{s} \in \mathbb{Z}^n, \ \underline{u} \in \mathbb{C}^n$ and $\underline{r} \in \Gamma$. It can be verified that $L(\underline{\beta},V_1)$ is indeed a $\widehat{\tau}_{0}$-module (see \cite[Section 8]{BE}).\\    
For $v \in V_{+}$, let us denote its image in $\widetilde{V}$ by $\overline{v}$ and define 
\begin{align*}
\widetilde{\phi} : V_{+} \longrightarrow L(\underline{\beta}, \widetilde{V})
\end{align*}
$\,\,\,\,\,\,\,\,\,\,\,\,\,\,\,\,\,\,\,\,\,\,\,\,\,\,\,\,\,\,\,\,\,\,\,\,\,\,\,\,\,\,\,\,\,\,\,\,\,\,\,\,\,\,\,\,\,\,\,\,\,\,\,\,\,\,\,\,\,\,\,\,\,\,\,\,\,\,\,\,\,\,\,\,\ v \longmapsto \ \overline{v} \otimes t^{\underline{k}},\ v \in V_{+}(\underline{k})$.\\
This map is clearly a non-zero $\widehat{\tau}_{0}$-module homomorphism and thus from the irreducibility of $V_{+}$, it follows that $V_{+} \cong \widetilde{\phi}(V_{+})$ is a $\widehat{\tau}_{0}$-submodule of $L(\underline{\beta}, \widetilde{V})$.\\
Recall that $\widetilde{V}= \bigoplus_{\overline{p} \in G} \widetilde{V}_{\overline{p}}$ and for every $\overline{p} \in G$, set 
\begin{align*}
L(\underline{\beta}, \widetilde{V})(\overline{p}) = \text{span} \{v \otimes t^{\underline{k} + \underline{r} + \underline{p}} \ | \ \underline{r} \in \Gamma, \ \underline{k} \in \mathbb{Z}^n, \ v \in \widetilde{V}_{\overline{k}} \}    
\end{align*}
which thereby forms a $\widehat{\tau}_{0}$-module under the usual action induced from (\ref{Action}).\\
The following results can be deduced similarly as in \cite{BE} by 
simply working with our Lie algebra $\widehat{\tau}_{0}$ instead of $L$.     
\bppsn\label{P7.3}
\
\begin{enumerate}
\item $V_{+} \cong L(\underline{\beta}, \widetilde{V})(\overline{0})$ as $\widehat{\tau}_{0}$-modules.
\item $\widetilde{V}$ is a $G$-graded-irreducible module over $\widetilde{L}$.
\item $\widetilde{V}$ is a completely reducible module over $\widetilde{L}$ and all its irreducible components are mutually isomorphic as $\big (I \ltimes L(\fma_{\overline{0}}, \sigma) \big )$-modules.        
\end{enumerate}
\eppsn

For any $d \in \mathbb{N}, \ \underline{k} \in \mathbb{Z}^n, \ x \in \fma_{\overline{k}}$ and $\underline{r_1}, \ldots, \underline{r_d}  \in \Gamma$, define $x(\underline{k}, \underline{r_1}, \ldots, \underline{r_d}) := \\ x \otimes t^{\underline{k}} - \sum_{i} x \otimes t^{\underline{k} + \underline{r_i}} + \sum_{i<j} x \otimes t^{\underline{k} + \underline{r_i} + \underline{r_j}} + \ldots + (-1)^d x \otimes t^{\underline{k} + \underline{r_1} + \underline{r_2} + \ldots + \underline{r_d}}$ and set $F_d$= span$\{x(\underline{k}, \underline{r_1}, \ldots, \underline{r_d}) \ | \ x \in \fma_{\overline{k}}, \ \underline{k} \in \mathbb{Z}^n, \ \underline{r_1}, \ldots, \underline{r_d} \in \Gamma \}$. We can easily check that the subspace $F_d$ is an ideal of $L(\fma, \sigma)$. \\
Recall that $I$ = span$\{I(\underline{u}, \underline{r}) \ | \ \underline{u} \in \mathbb{C}^n, \ \underline{r} \in \Gamma \}$ and define for any $d \in \mathbb{N}$, 
\begin{align*}
 I_d(\underline{u}, \underline{r}, \underline{s_1}, \ldots,  \underline{s_d}):=I(\underline{u}, \underline{r}) - \sum_{i} I(\underline{u}, \underline{r}+ \underline{s_i})+ \sum_{i<j} I(\underline{u}, \underline{r}+ \underline{s_i}+ \underline{s_j})\\  + \ldots + (-1)^{d} I(\underline{u}, \underline{r}+ \underline{s_1} + \underline{s_2} + \underline{s_d}) 
\end{align*}
where $\underline{u} \in \mathbb{C}^n$ and $\underline{r}, \ \underline{s_1}, \ldots, \underline{s_d} \in \Gamma$.\\  
Let $I_d$ = span$\{I_d(\underline{u}, \underline{r}, \underline{s_1}, \ldots, \underline{s_d}) \ | \ \underline{u} \in \mathbb{C}^n, \ \underline{r}, \ \underline{s_1}, \ldots, \underline{s_d} \in \Gamma \}$. Then it can be verified that $I_d$ is an ideal in $I$. Further consider the map 
\begin{align}\label{pi}
\pi : I/I_2 \longrightarrow \mathfrak{gl}_n(\mathbb{C})
\end{align}
$\,\,\,\,\,\,\,\,\,\,\,\,\,\,\,\,\,\,\,\,\,\,\,\,\,\,\,\,\,\,\,\,\,\,\,\,\,\,\,\,\,\,\,\,\,\,\,\,\,\,\,\,\ I(e_i, m_je_j) + I_2  \longmapsto m_jE_{ji} \ , \ 1 \leqslant i,j \leqslant n$\\
and extend $\pi$ linearly. Here $\{e_i\}_{i=1}^{n}$ denotes  the standard basis of $\mathbb{C}^n$. It is not very hard to deduce from \cite[Proposition 4.1]{E2} that $I/I_2 \cong \mathfrak{gl}_n(\mathbb{C})$ as Lie algebras and every finite-dimensional irreducible $I$-module eventually factors through $\mathfrak{gl}_n(\mathbb{C})$ via this map $\pi$.\\
Again if we consider the map 
$ev_{\underline{1}} : L(\fma, \sigma) \longrightarrow \fma$

$\,\,\,\,\,\,\,\,\,\,\,\,\,\,\,\,\,\,\,\,\,\,\,\,\,\,\,\,\,\,\,\,\,\,\,\,\,\,\,\,\,\,\,\,\,\,\,\,\,\,\,\,\,\,\,\,\,\,\,\,\,\,\,\,\,\,\,\,\,\,\,\,\,\,\,\,\,\,\,\,\,\ x \otimes t^{\underline{k} + \underline{r}} \longmapsto \ x ,\ x \in \fma_{\overline{k}} \ (\underline{k} \in \mathbb{Z}^n, \underline{r} \in \Gamma)$, \\ 
then by means of some elementary arguments involving linear independence, it is not difficult to show that $Ker$($ev_{\underline{1}}$)=$F_1$. But since this map is onto, we clearly have $L(\fma, \sigma)/F_1 \cong \fma$ as Lie algebras.     
 
\blmma \label{L7.4}  
Let $\underline{k} \in \mathbb{Z}^n, \ x \in \fma_{\overline{k}} \ , \ \underline{s_1}, \ldots, \underline{s_d} \in \Gamma$ and $d \in \mathbb{N}$. If the element $x(\underline{k}, \underline{s_1}, \ldots, \underline{s_d}) \in F_d$ acts by the scalar $\lambda(\underline{k}, \underline{s_1}, \ldots, \underline{s_d})$ on a finite-dimensional representation $(\rho^{\prime},V^{\prime})$ of $\widetilde{L}$, then $\lambda(\underline{k}, \underline{s_1}, \ldots, \underline{s_d})=0$.   
\elmma
\begin{proof}
Let $\underline{r} \in \Gamma$ be arbitrary. Now we have
\begin{align*}
[I_d(\underline{u}, \underline{r}, \underline{s_1}, \ldots,  \underline{s_d}), x \otimes t^{\underline{k}}] = (\underline{u}, \underline{k})x(\underline{k} + \underline{r},\underline{s_1}, \ldots, \underline{s_d}) \ \forall \ \underline{u} \in \mathbb{C}^n.         
\end{align*}
But as $V^{\prime}$ is finite-dimensional, we can simply take the trace operator after applying the map $\rho^{\prime}$ on both sides to finally obtain 
\begin{align*}
\lambda(\underline{k} + \underline{r},\underline{s_1}, \ldots, \underline{s_d})=0 \ \forall \ \underline{k} \in \mathbb{Z}^n \setminus \{\underline{0}\}, \ \underline{r} \in \Gamma.   
\end{align*} 
This immediately gives us the desired conclusion.              
\end{proof}

We already know from our Proposition \ref{P7.3} that there exist irreducible $\widetilde{L}$-modules $\widetilde{M}_{1}, \ldots, \widetilde{M}_{N}$ such that 
$\widetilde{V}= \bigoplus_{\nu=1}^{N} \widetilde{M}_{\nu}$ for some $N \in \mathbb{N}$ with $\widetilde{M}_{i} \cong \widetilde{M}_{j}$ as $\big (I \ltimes L(\fma_{\overline{0}}, \sigma) \big )$-modules $\forall \ 1 \leqslant i,j \leqslant N$. It is also clear that we have $L(\fma_{\overline{0}}, \sigma) = \mathfrak{h}_{\overline{0}} \otimes A(\underline{m})$.                           

\bthm
For each $1 \leqslant \nu \leqslant N, \ \widetilde{M}_{\nu}$ is a finite-dimensional irreducible module for the direct sum $\mathfrak{gl}_n(\mathbb{C}) \oplus \fma$. Thus $\widetilde{V}$ is the unique completely reducible finite-dimensional $G$-graded-irreducible quotient of $V_{+}$ over $\widetilde{L}$.  
\ethm
\begin{proof}
Using our Lemma \ref{L7.4}, we can give a similar argument as presented in \cite[Theorem 9.4]{BE} to show that the Lie algebra $I_2 \oplus F_1$ actually acts trivially on every $\widetilde{M}_{\nu}$. The theorem now directly follows from (\ref{pi}). The final assertion is clear from the definitions of $F_1$ and the $\widetilde{L}$-submodule $W$ of $V_{+}$.             
\end{proof}
\brmk
The second author would like to take this opportunity to correct some results in \cite{BE}. The result stated in \cite[Proposition 8.12(1)]{BE} is not correct as the eigenspaces need not be isomorphic as $\mathfrak{g}$-modules, rather they are mutually isomorphic only as $\mathfrak{g}_{\overline{0}}$-modules where $\mathfrak{g}$ is a graded Lie algebra. This is simply because the map defined at the beginning of Page-99 in \cite{BE} is a $\mathfrak{g}_{\overline{0}}$-module homomorphism and not necessarily a  $\mathfrak{g}$-module homomorphism as claimed by the authors. As a result, Theorem 8.3 in \cite{BE} is also affected and we ultimately obtain $\widetilde{V}$ to be a completely reducible $\widetilde{L}$-module with its irreducible components being isomorphic as $\big (I \ltimes L(\fma_{\overline{0}}, \sigma) \big )$-modules. Finally although \cite[Lemma 9.5]{BE} is correct, the proof provided there is not valid in this context due to the reasons stated above and thereby cannot be employed to prove \cite[Theorem 9.4]{BE}. Nonetheless, the trace argument given in our Lemma \ref{L7.4} can be utilized to resolve this problem and complete the proof of Theorem 9.4 in \cite{BE}. We would like to also point out that \cite[Theorem 8.3]{BE} is again recalled in \cite{BES} and \cite{E3}, but not effectively used in both these papers.                                           

\ermk
For each $1 \leqslant \nu \leqslant N$, consider any surjective Lie algebra homomorphism 
\begin{align*}
\phi_{\nu} : \widetilde{L} \longrightarrow \mathfrak{gl}_n(\mathbb{C}) \oplus \fma
\end{align*}
satisfying $Ker \phi_{\nu}=F_1 \oplus I_2$ and which maps $L(\fma, \sigma)$ onto $\fma$ and $I$ onto $\mathfrak{gl}_n(\mathbb{C})$. Set $\phi_{\nu}|_{I}= \phi_{\nu}^{\prime}$ and $\phi_{\nu}|_{L(\fma, \sigma)}= \phi_{\nu}^{\prime \prime}$. Using this identification  map $\phi_{\nu}$ and then applying \cite[Lemma 2.7]{H}, we can thereby find finite-dimensional irreducible representations, say $\widetilde{V}_{1}^{\nu}$ and $\widetilde{W}_{2}^{\nu}$ of $\mathfrak{gl}_n(\mathbb{C})$ and $\fma$ respectively so that we have an isomorphism of irreducible $\big (\mathfrak{gl}_n(\mathbb{C}) \oplus \fma \big)$-modules given by
\begin{align*}
\widetilde{M}_{\nu} \cong \widetilde{V}_{1}^{\nu} \otimes \widetilde{W}_{2}^{\nu} \ \forall \ 1 \leqslant \nu \leqslant N.
\end{align*}
Note that due to (\ref{pi}), it suffices to only consider the same identification map $\phi_{\nu}^{\prime}= \pi$ for defining a $I$-module structure on all the $\widetilde{V}_{1}^{\nu}$'s. Consequently all the $\mathfrak{gl}_n(\mathbb{C})$-modules $\widetilde{V}_{1}^{\nu}$ remain isomorphic to each other which permits us to take $\widetilde{V}_{1}^{\nu} \cong \widetilde{V}_{1}$ as $\mathfrak{gl}_{n}(\mathbb{C})$-modules for all $1 \leqslant \nu \leqslant N$. Also $L(\fma, \sigma)$ sits diagonally inside $L(\fma, \sigma)^{\oplus{N}}$ and $\bigoplus_{\nu=1}^{N} \phi_{\nu}^{\prime \prime}$ is a Lie algebra homomorphism that maps $L(\fma, \sigma)^{\oplus{N}}$ onto $\fma^{\oplus{N}}$. This immediately induces an $L(\fma, \sigma)$-module structure on  $\widetilde{W}_{2}:= \sum_{\nu=1}^{N}\widetilde{W}_{2}^{\nu}$. Moreover since all the $\widetilde{W}_{2}^{\nu}$'s are irreducible, we can assume (without any loss of generality) that the sum given in $\widetilde{W}_{2}$ is $direct$. We can now invoke Proposition \ref{P7.3} to infer that $\widetilde{W}_{2}$ is a $G$-graded-irreducible $L(\fma, \sigma)$-module where we are allowed to take $\widetilde{V}_{1}$ to be $zero$-$graded$ as $I$ lies in the zeroth graded component of $\widetilde{L}$.\\
We now define a $\widehat{\tau}_0$-module structure on $\widetilde{V}_{1} \otimes \widetilde{W}_{2} \otimes A$ by setting
\begin{align*}
(x \otimes t^{\underline{l}}). \big (\sum_{\nu=1}^{N}v_1 \otimes w_{2}^{\nu} \otimes t^{\underline{k}} \big )= \sum_{\nu=1}^{N}v_1 \otimes \big ( \phi_{\nu}^{\prime \prime}(x \otimes t^{\underline{l}}).w_{2}^{\nu} \big ) \otimes t^{\underline{k} + \underline{l}} \ ;\\
t^{\underline{r}}d_i.\big (\sum_{\nu=1}^{N}v_1 \otimes w_{2}^{\nu} \otimes t^{\underline{k}} \big )=(k_i + \beta_i) \big (\sum_{\nu=1}^{N}v_1 \otimes w_{2}^{\nu} \otimes t^{\underline{k} + \underline{r}} \big )\\
+ \ \sum_{\nu=1}^{N} \sum_{j=1}^{n} \big ((r_jm_jE_{ji}).v_1 \big ) \otimes w_{2}^{\nu} \otimes t^{\underline{k} + \underline{r}}
\end{align*}
for all $x \in \fma_{\overline{l}}, \ \underline{k}, \ \underline{l} \in \mathbb{Z}^n, \ \underline{r} = \sum_{j=1}^{n}r_jm_je_j \in \Gamma, \ v_1 \in \widetilde{V}_{1}, \ w_{2}^{\nu} \in \widetilde{W}_{2}^{\nu}$ and $1 \leqslant i \leqslant n$. Then it can be trivially verified that $\widetilde{V}_{1} \otimes \widetilde{W}_{2} \otimes A$ is a  $\widehat{\tau}_0$-module. Let $\widetilde{W}_{2}=  \bigoplus_{\overline{k} \in G} \widetilde{W}_{2, \overline{k}}$ and define $T^{\prime}:= \bigoplus_{\underline{k} \in \mathbb{Z}^n} \widetilde{V}_{1} \otimes \widetilde{W}_{2, \overline{k}} \otimes \mathbb{C}t^{\underline{k}}$ which is a submodule of $\widetilde{V}_{1} \otimes \widetilde{W}_{2} \otimes A$. One can now check that $T^{\prime}$ is in fact an irreducible module over $\widehat{\tau}_0$ and henceforth directly define a natural non-zero $\widehat{\tau}_0$-module homomorphism between $L(\underline{\beta}, \widetilde{V})(\overline{0})$ and $T^{\prime}$. This subsequently becomes an isomorphism due to the irreducibility of both $L(\underline{\beta}, \widetilde{V})(\overline{0})$ and $T^{\prime}$.  Finally using Proposition \ref{P7.3}, it is trivial to see that $V_{+} \cong T^{\prime}$ as $\widehat{\tau}_0$-modules.\\
Consider the following triangular decomposition of $LT$ (see (\ref{LT}) and  (\ref{Grading})) given by 
\begin{align*}
LT = \widetilde{\tau}_{-} \oplus L(\fma, \sigma) \oplus \widetilde{\tau}_{+}.
\end{align*}
Let $\widetilde{V}_{2}$ stand for the unique $G$-graded-irreducible quotient of the induced module for $\widetilde{W}_{2}$ with respect to this triangular decomposition. Also for each $1 \leqslant \nu \leqslant N$, let us denote the unique irreducible quotient of the induced representation of $\widetilde{W}_{2}^{\nu}$ for the same triangular decomposition by $(\widetilde{V}_{2}^{\nu}, \widetilde{\rho}_{\nu})$. Now
we already know that $\widetilde{W}_{2}= \bigoplus_{\nu=1}^{N} \widetilde{W}_{2}^{\nu}$ as $L(\fma, \sigma)$-modules. Then by analogous arguments presented in \cite[Lemma 3.2]{E3}, we can easily infer that $\widetilde{V}_{2}= \bigoplus_{\nu=1}^{N}\widetilde{V}_{2}^{\nu}$ and thus  $\bigoplus_{\nu=1}^{N}\widetilde{V}_{2}^{\nu}$ is a $G$-graded-irreducible $LT$-module.

Again define a $\widehat{\tau}$-module structure on $\widetilde{V_1} \otimes \widetilde{V_2} \otimes A$ by setting
\begin{align*}
(x \otimes t^{\underline{l}}). \big (\sum_{\nu=1}^{N}v_1 \otimes v_{2}^{\nu} \otimes t^{\underline{k}} \big )= \sum_{\nu=1}^{N}v_1 \otimes \big ( \widetilde{\rho}_{\nu}(x \otimes t^{\underline{l}})v_{2}^{\nu} \big ) \otimes t^{\underline{k} + \underline{l}} \ ;\\
t^{\underline{r}}d_i.\big (\sum_{\nu=1}^{N}v_1 \otimes v_{2}^{\nu} \otimes t^{\underline{k}} \big )=(k_i + \beta_i) \big (\sum_{\nu=1}^{N}v_1 \otimes v_{2}^{\nu} \otimes t^{\underline{k} + \underline{r}} \big )\\
+ \ \sum_{\nu=1}^{N} \sum_{j=1}^{n} \big ((r_jm_jE_{ji}).v_1 \big ) \otimes v_{2}^{\nu} \otimes t^{\underline{k} + \underline{r}}
\end{align*}
for all $x \in \mathfrak{g}_{\overline{l}}, \ \underline{k}, \ \underline{l} \in \mathbb{Z}^n, \ \underline{r} = \sum_{j=1}^{n}r_jm_je_j \in \Gamma, \ v_1 \in \widetilde{V}_{1}, \ v_{2}^{\nu} \in \widetilde{V}_{2}^{\nu}$ and $1 \leqslant i \leqslant n$. One can readily check that $\widetilde{V}_{1} \otimes \widetilde{V}_{2} \otimes A$ is a  $\widehat{\tau}$-module. Put $\widetilde{V}_{2}=  \bigoplus_{\overline{k} \in G} \widetilde{V}_{2, \overline{k}}$ and define $V^{\prime}:= \bigoplus_{\underline{k} \in \mathbb{Z}^n} \widetilde{V}_{1} \otimes \widetilde{V}_{2, \overline{k}} \otimes \mathbb{C}t^{\underline{k}}$ which is again a submodule of $\widetilde{V}_{1} \otimes \widetilde{V}_{2} \otimes A$.
\bppsn\label{P7.7}
$V^{\prime}$ is an irreducible module over $\widehat{\tau}$.
\eppsn
\begin{proof}
The proof is similar to Proposition \ref{P6.1}. 
\end{proof}	

\bthm \label{T7.13}
$V \cong V^{\prime}$ as $\widehat{\tau}$-modules.
\ethm
\begin{proof}
It is immediate that the highest weight space of $V^{\prime}$ is given by $T^{\prime}$. The result now directly follows from Proposition \ref{P7.7} as
$V_{+} \cong T^{\prime}$. 
\end{proof}

\bcrlre\label{C7.11} 
There exists a quadruplet $(\psi, c, \lambda, \underline{\beta}) \in P_{\mathfrak{sl}_n}^{+} \times \mathbb{C} \times (P_{\mathfrak{g}}^{+})^{\times} \times \mathbb{C}^n$ such that $V \cong V(\psi, c, \lambda, \underline{\beta})$ as $\widetilde{\tau}$-modules.     
\ecrlre
\begin{proof}
From our construction of $V^{\prime}$ and Theorem \ref{T7.13}, it is evident that the integrability of $V$ also forces each $\widetilde{V}_{2, \overline{k}}$ to be integrable over the finite-dimensional simple Lie algebra $\mathfrak{g}_{\overline{0}}$. Moreover as $\widetilde{W}_{2}$ is finite-dimensional, we obtain that $\widetilde{V}_{2, \overline{k}}$ has finite-dimensional weight spaces with respect to $\mathfrak{h}_{\overline{0}}$. Subsequently by \cite[Lemma 3.5]{S}, it readily follows that $\widetilde{V}_{2, \overline{k}}$ must be finite-dimensional which shows that $\widetilde{V}_{2}$ is also finite-dimensional as $G$ is a finite group. In particular, $\widetilde{V}_{2}^{\nu}$ is a finite-dimensional irreducible $LT$-module for each $1 \leqslant \nu \leqslant N$. Then by \cite[Corollary 4.4]{L}, we can directly conclude that there exists a finite-dimensional irreducible representation $(V_{\nu},f_{\nu})$ of $\mathfrak{g}$ and $\underline{a_{\nu}} \in (\mathbb{C}^{\times})^n$ such that $\widetilde{V}_{2}^{\nu} \cong ev_{\underline{a_{\nu}}} V_{\nu}$ as $LT$-modules for all $1 \leqslant \nu \leqslant N$. This eventually implies that it suffices to only consider $\widetilde{\rho}_{\nu}$ to be the composition of $ev_{\underline{a_{\nu}}}$ and $f_{\nu}$ which consequently gives $\phi_{\nu}^{\prime \prime} = ev_{\underline{a_{\nu}}}$. But again as $Ker\phi_{\nu}^{\prime \prime}=F_1$, it is clear that $\underline{a_{\nu}} \in \mathcal{U}(\underline{m}) \ \forall \ 1 \leqslant \nu \leqslant N$. Moreover since we know that $\widetilde{V}_{2}$ is a $G$-graded-irreducible module over $LT$, we can infer that the twisted $\mathfrak{g}$-module under the action 
\begin{align*}
x. \big (\sum_{\nu=1}^{N} v_{2}^{\nu} \big)= \sum_{\nu=1}^{N} \underline{a_{\nu}}^{\underline{l}}(x.v_{2}^{\nu}) \ \forall \ x \in \mathfrak{g}_{\overline{l}}, \ v_{2}^{\nu} \in V_{2}^{\nu}
\end{align*}
is also $G$-graded-irreducible. The corollary then follows from Theorem \ref{T7.13} and Remark \ref{R6.2}. Lastly observe that since $LT$ does not act trivially on $V$, we can $never$ have $\lambda=0$.   
\end{proof}

\section{Isomorphism Classes of Irreducible Modules}\label{Class}  
In this section, our primary goal is to determine the isomorphism classes of the irreducible integrable modules that we have classified in Section \ref{Classification}.
          
We first note that in our construction of $V(\psi, c, \lambda, \underline{\beta})$ in Section \ref{Realization}, the vector space $V_1 \otimes V_2$ is endowed with a $(I \ltimes LT)$-module structure given by
\begin{align*}
(x \otimes t^{\underline{l}}). (v_1 \otimes v_2)= v_1 \otimes (x.v_2), \\
I(\underline{u},\underline{r}).(v_1 \otimes v_2)= \sum_{i,j=1}^{n} \big ((u_ir_jm_jE_{ji}).v_1 \big ) \otimes v_2
\end{align*}
for all $x \in \mathfrak{g}_{\overline{l}}, \ \underline{l} \in \mathbb{Z}^n, \ \underline{u} \in \mathbb{C}^n, \ \underline{r} = \sum_{j=1}^{n}r_jm_je_j \in \Gamma$ and $v_1 \in V_1, \ v_2 \in V_2$.              
\bppsn \label{P8.2}
Let $\bigoplus_{\underline{k} \in \mathbb{Z}^n}V_{1} \otimes V_{2, \overline{k}} \otimes \mathbb{C}t^{\underline{k}} \cong \bigoplus_{\underline{k} \in \mathbb{Z}^n}V_{1}^{\prime} \otimes V_{2, \overline{k}}^{\prime} \otimes \mathbb{C}t^{\underline{k}}$ as $\widetilde{\tau}$-modules where $V_1, V_{1}^{\prime}$ are finite-dimensional irreducible $I$-modules and $V_2,  V_{2}^{\prime}$ are finite-dimensional $G$-graded-irreducible $LT$-modules used in
the construction of $V(\psi, c, \lambda, \underline{\beta})$ in Section \ref{Realization}. Then
\begin{enumerate}
\item $V_1 \otimes V_2 \cong V_{1}^{\prime} \otimes V_{2}^{\prime}$ as  $(I \ltimes LT)$-modules.
\item $V_1 \cong V_{1}^{\prime}$ as $\mathfrak{gl}_{n}(\mathbb{C})$-modules and $V_2 \cong V_{2}^{\prime}$ as $G$-graded-irreducible $\mathfrak{g}$-modules up to a possible $shift$ $of$ $grading$. 
\end{enumerate}
\eppsn
\begin{proof}

(1) Let $W_2$ and $W_{2}^{\prime}$ be the highest weight spaces of the $G$-graded-irreducible $LT$-modules $V_2$ and $V_{2}^{\prime}$ respectively.  We claim that both the highest weight spaces are also $G$-graded-irreducible modules over $L(\fma, \sigma)$. It suffices to prove our claim only for $W_2$. As $V_2$ is $G$-graded, we have 
$w= \sum_{\overline{k} \in G}w_{\overline{k}}$ for any $w \in W_2$. Applying an element of the positive root space $\widetilde{\tau}_{+}$ (see (\ref{LT})) on $w$, it is easy to see that $w_{\overline{k}} \in W_2$ for all $\overline{k} \in G$ which in turn shows that $W_2$ is also $G$-graded. The claim now directly follows from the $G$-graded-irreducibility of $V_2$ and the PBW theorem. Now the isomorphism between the irreducible $\widetilde{\tau}$-modules again gives rise to an isomorphism of the highest weight spaces, say $V_{+}$ and $V_{+}^{\prime}$ over $\widetilde{\tau}_0$   
\begin{align*}
F : \bigoplus_{\underline{k} \in \mathbb{Z}^n}V_1 \otimes W_{2, \overline{k}} \otimes \mathbb{C}t^{\underline{k}} = V_{+} \longrightarrow V_{+}^{\prime} = \bigoplus_{\underline{k} \in \mathbb{Z}^n}V_{1}^{\prime} \otimes W_{2, \overline{k}}^{\prime} \otimes \mathbb{C}t^{\underline{k}} \ .                 
\end{align*}
Moreover as $W_2$ and $W_{2}^{\prime}$ are finite-dimensional $G$-graded-irreducible modules over $L(\fma, \sigma)$ and $[\mathfrak{h}_{\overline{0}}, \fma]=(0)$, it follows that there exists a unique $\overline{\lambda} \in \mathfrak{h}_{\overline{0}}^{*}$ such that $P(V_{+})= \{ \overline{\lambda} + \delta_{\underline{k} + \underline{\beta}} \ | \ \underline{k} \in \mathbb{Z}^n \}$ and $P(V_{+}^{\prime})=\{ \overline{\lambda} + \delta_{\underline{k} + \underline{\beta^{\prime}}} \ | \ \underline{k} \in \mathbb{Z}^n \}$ for some $\underline{\beta}, \ \underline{\beta^{\prime}} \in \mathbb{C}^n$ with $\underline{\beta}-\underline{\beta^{\prime}} \in \mathbb{Z}^n$. From our construction, it is also evident that $L(\fma, \sigma)$ does not act trivially on both $W_2$ and $W_{2}^{\prime}$ and thus we must have $\overline{\lambda} \neq0$. For the sake of notational convenience, we shall assume, without loss of generality, that $\underline{\beta}= \underline{0}= \underline{\beta^{\prime}}$. Again as $F$ is an isomorphism of $\widetilde{\tau}_0$-modules, the $\underline{k}$-weight vectors in $V_{+}$ should map to $\underline{k}$-weight vectors in $V_{+}^{\prime}$ and therefore for some $l \in \mathbb{N}$, we must have 
\begin{align*}
F (v_1 \otimes w_2 \otimes t^{\underline{k}})=
\sum_{i=1}^{l}v_{1,i}^{\prime} \otimes w_{2,i}^{\prime} \otimes t^{\underline{k}} \ , \ v_{1,i}^{\prime} \in V_{1}^{\prime}, \ w_{2,i}^{\prime} \in W_{2, \overline{k}}^{\prime} \ , \ \underline{k} \in \mathbb{Z}^n.                            
\end{align*}	
Define $ \phi : V_1 \otimes W_2 \longrightarrow V_{1}^{\prime} \otimes W_{2}^{\prime}$ by setting
\begin{align*}
 v_1 \otimes w_2 \mapsto \sum_{i=1}^{l} v_{1,i}^{\prime} \otimes w_{2,i}^{\prime} , \ v_1 \in V_1, \ v_{1,i}^{\prime} \in V_{1}^{\prime} \  \text{and} \ w_2 \in W_{2, \overline{k}} \ , \ w_{2,i}^{\prime} \in W_{2, \overline{k}}^{\prime}.                                    
\end{align*}
Now using the fact that $\overline{\lambda} \neq 0$ and the relation\\   
$F \big ((h \otimes t^{\underline{r}}).(v_1 \otimes w_2 \otimes t^{\underline{k}}) \big ) = (h \otimes t^{\underline{r}}).F \big (v_1 \otimes w_2 \otimes t^{\underline{k}} \big ) \ \forall \ h \in \mathfrak{h}_{\overline{0}}$ and $\underline{r} \in \Gamma$, \\ it is trivial to verify that $\phi$ is well-defined. Furthermore, it is also clear that $\phi$ is in fact an isomorphism of $(I \ltimes L(\fma, \sigma))$-modules. But since the $G$-graded-irreducible modules $V_1 \otimes W_2$ and $V_{1}^{\prime} \otimes W_{2}^{\prime}$ are clearly the highest weight spaces of the $G$-graded-irreducible modules $V_{1} \otimes V_{2}$ and $V_{1}^{\prime} \otimes V_{2}^{\prime}$ respectively, this finally yields 
$V_1 \otimes V_{2} \cong V_{1}^{\prime} \otimes V_{2}^{\prime}$ as $(I \ltimes LT)$-modules. Note that this isomorphism is also 
a homomorphism of $G$-graded modules over $(I \ltimes LT)$.\\ 
(2) We already have $V_1 \otimes W_2 \cong V_{1}^{\prime} \otimes W_{2}^{\prime}$ as $(I \ltimes L(\fma, \sigma))$-modules. Let us refer to this isomorphism as $\psi$. Now choose bases $\{w_1, \ldots, w_r \}$ and $\{w_{1}^{\prime}, \ldots, w_{s}^{\prime} \}$ of $W_2$ and $W_{2}^{\prime}$ respectively. Then we obtain an isomorphism of $I$-modules 
\begin{align*}
\bigoplus_{i=1}^{r} V_1 \otimes \mathbb{C}w_i \cong \bigoplus_{j=1}^{s} V_{1}^{\prime} \otimes \mathbb{C}w_{j}^{\prime}
\end{align*}
whence we get $r=s$ and $V_1 \cong V_{1}^{\prime}$ as $\mathfrak{gl}_{n}(\mathbb{C})$-modules. Again let us pick a basis  
$\{v_1, \ldots, v_k \}$ of $V_1$ and set $\psi_i= \psi|_{\mathbb{C}v_i \otimes W_2}$ for each $1 \leqslant i \leqslant k$. Then there exists some
$1 \leqslant i \leqslant k$ such that $\psi_i \neq 0$. Composing $\psi_i$ with the obvious projection map from $V_{1}^{\prime} \otimes W_{2}^{\prime}$ to $W_{2}^{\prime}$, we get a non-zero graded-module homomorphism, say $\epsilon$, between the $G$-graded-irreducible $L(\fma, \sigma)$-modules $\mathbb{C}v_i \otimes W_2$ and $W_{2}^{\prime}$ . This readily implies that $\epsilon$ must be injective. But as dim$W_2$ = dim$W_{2}^{\prime}$ , it follows that $\epsilon$ is an $L(\fma, \sigma)$-module isomorphism. Thus we have $V_2 \cong V_{2}^{\prime}$ as $LT$-modules which in turn gives us the desired result from Proposition \ref{P5.1}.    
\end{proof}

\bthm\label{T8.4}
If $V(\psi,c,\lambda, \underline{\beta}) \cong V(\psi^{\prime},c^{\prime},\lambda^{\prime}, \underline{\beta^{\prime}})$ as $\widetilde{\tau}$-modules, then
\begin{enumerate}
\item $\psi= \psi^{\prime}$ and $c=c^{\prime}$;
\item $\lambda^{\prime} \in \widehat{G}\lambda$ which in particular shows that $\widehat{G}\lambda=\widehat{G}\lambda^{\prime}$ and $S(\lambda)=S(\lambda^{\prime})$;
\item $\underline{\beta}- \underline{\beta^{\prime}} \in \mathbb{Z}^n$.    
\end{enumerate}
\ethm
\begin{proof}
This follows directly from Subsection \ref{SS5.2} and Proposition \ref{P8.2}.         
\end{proof}

\section{The Final Theorem}\label{Final}          
Let us define an equivalence relation '$\sim$' on $\mathbb{C}^n$ by declaring $\underline{\beta} \sim \underline{\beta^{\prime}} \iff \underline{\beta} - \underline{\beta^{\prime}} \in \mathbb{Z}^n$. Also define another equivalence relation '$\sim$' on $(P_{\mathfrak{g}}^{+})^{\times}$ by setting $\lambda \sim \lambda^{\prime} \iff \widehat{G}\lambda=\widehat{G}\lambda^{\prime}$. We shall denote the equivalence class under both these equivalence relations by $[ \  . \ ]$. We are now ready to state the main theorem of our paper which is a direct consequence of Proposition \ref{P6.1}, Remark \ref{R6.2}, Corollary \ref{C7.11} and Theorem \ref{T8.4}.                                            
\bthm
If $\mathcal{I}_{fin}^{0}$ denotes the set of all level zero irreducible integrable modules over $\widetilde{\tau}$ (upto isomorphism) having finite-dimensional weight spaces with respect to $\widetilde{\mathfrak{h}}$ and endowed with a non-trivial $LT$-action, then we have a well-defined bijective map
\begin{align*} 
\kappa \ : \ P_{\mathfrak{sl}_n}^{+} \times \mathbb{C} \times (P_{\mathfrak{g}}^{+})^{\times}/\sim \times \ \mathbb{C}^n/\sim \ \longrightarrow \ \mathcal{I}_{fin}^{0}    
\end{align*}
\,\,\,\,\,\,\,\,\,\,\,\,\,\,\,\,\,\,\,\,\,\,\,\,\,\,\,\,\,\,\,\,\,\,\,\,\,\,\,\,\,\,\,\,\,\,\,\,\,\,\,\,\,\,\,\,\,\,\,\,\ $(\psi, c, [\lambda], [\underline{\beta}] ) \,\,\,\,\,\,\,\,\,\,\,\ \mapsto \,\,\,\,\,\,\,\,\,\,\,\,\,\,\,\,\,\,\ V(\psi, c, [\lambda], [\underline{\beta}])$.                          
\ethm 

\brmk
Our Corollary \ref{C7.11} also recovers the classification result given in \cite[Theorem 4.3]{EJ} for the untwisted case. In fact, we have further deduced that any such level zero irreducible integrable module is uniquely determined by the parameters $(\psi, c, \lambda, [\underline{\beta}]) \in P_{\mathfrak{sl}_n}^{+} \times \mathbb{C} \times (P_{\mathfrak{g}}^{+})^{\times} \times \mathbb{C}^n/\sim$. Note that since  $S(\lambda)=1$ and $|\widehat{G}\lambda|=1$, these parameters (which are dependent on $\lambda$) do not make any appearance in the description of these irreducible integrable modules in the untwisted case.                                          
\ermk

\noindent $\bf{Acknowledgements.}$ The first author would like to thank Prof. Mikhail Kochetov for some helpful discussions. The authors are also grateful to the anonymous referee for useful comments and suggestions.

\end{document}